\documentclass [11pt]{article}

\usepackage{graphicx}

\newtheorem {theo} {\bf Theorem} [section]
\newtheorem {prop} [theo] {\bf Proposition}
\newtheorem {cory} [theo]{\bf Corollary}
\newtheorem {lem} [theo] {\bf Lemma}

\newtheorem {obs}  [theo]{\bf Remark}

\newcommand{\CaixaPreta}{\rule{2mm}{2mm}}
\newcommand{\qed}{\hfill\caixapreta \vspace{5mm}}

\newtheorem {REMCURSIVA} [theo] {\bf Remark}

\newcommand{\be}{\begin{eqnarray}}
\newcommand{\ee}{\end{eqnarray}}
\newcommand{\benn}{\begin{eqnarray*}}
\newcommand{\eenn}{\end{eqnarray*}}
\newcommand{\bse}{\begin{equation}}
\newcommand{\ese}{\end{equation}}
\newcommand{\bsenn}{\begin{displaymath}}
\newcommand{\esenn}{\end{displaymath}}

\newcommand{\N}{\mbox{I${\!}$N}}
\newcommand{\R}{\mbox{I${\!}$R}}

\newcommand{\coes}{\c c\~oes }

\newcommand{\rv}{{\bf \mbox {r}}}
\newcommand{\vv}{{\bf \mbox {v}}}

\newcommand{\nn}{{\bf \mbox {n}}}
\newcommand{\sv}{{\bf \mbox {s}}}

\title{On the Fixed Homogeneous Circle Problem}
 \author{C. Azev\^edo, H. Cabral and P. Ontaneda\footnote{The first author was
 partially supported by a CNPq/FACEPE research grant. The third author was
 partially supported by a research grant from CNPq, Brazil.} }
\date{}

\usepackage{wrapfig}
\begin{document}
\maketitle


\begin{abstract} We give some results about the dynamics of a particle moving in
Euclidean
three-space under the influence of the gravitational force induced by a fixed
homogeneous circle. Our main results concern (1) singularities and (2) the dynamics
in the plane that contains the circle. The study presented here is purely analytic.
\end{abstract}

In this paper we present a study about the movement of a particle
in Euclidean three-space $\R^{3}$ on which the only acting force is the
gravitational force induced by a fixed homogeneous circle.
Even though this  problem seems quite natural, we could not find in the literature
any reference concerning the dynamics of it. All we could find were a few different ways
of expressing the potential function. Essentially all these expressions had
already appeared in Poincare's Th\'eorie du Potentiel Newtonien
\cite{Po}, published first in 1899. Probably the most well-known expressions are
the one expressed in terms of elliptic integrals of the first kind
and the one using the arithmetic-geometric mean given by Gauss.
We believe this is the first study presented about the dynamics of
this problem.\\

Before we present our first result we need some notation and a couple of definitions.
We want to study the movement in $\R^{3}$ of a particle $P$
under the influence of the gravitational force
induced by a fixed homogeneous circle
$\cal C$.
Denote by ${\mbox{\rv}}=(x,y,z) \in \R^{3}-\cal C$ the position of the particle  $P$.
Also $\dot {\mbox{\rv}}=(\dot x,\dot y,\dot z)$ denotes its velocity.
According to Newton's Law the movement of $P$
obeys the following second order differential equation:

{\footnotesize\begin{equation}
\stackrel{..}{\mbox{\rv}} \,=\,-\nabla V({\mbox{\rv}}) \label{00}
\end{equation}}

\noindent where  $V$ denotes the potential energy
induced by $\,\cal C$. The expression of $V$ is given by
$V(\mbox{\rv})=-\int_{{\cal
C}}\frac{\lambda\,\,du}{|\!|{\mbox{\rv}}-u|\!|}$,
where  $\lambda$ is the constant mass density of the
circle $\cal C$. We say  that  a solution ${\mbox{\rv}}(t)$,
defined in the maximal interval  $(a,b)$, has a singularity
at $b$ (or at $a$), if $b<+\infty$ (if $a>-\infty$).
For $p\in \R^3$, let $dist(p,{\cal{C}})$ denote the Euclidean distance,
$inf_{x\in\cal C} \|p-x\|,$ from $p$ to $\cal C$.
We can now ask: if $lim_{t\rightarrow b^{-}}
dist({\mbox{\rv}}(t),{\cal{C}})=0,\,\,b<+\infty,$ does
$\,{\mbox{\rv}}(t)$  approach a well defined point in $\cal C$ ? A
priori, ${\mbox{\rv}}(t)$ could approach the circle without
getting closer to any specific point in the circle.
A singularity at $c$ is called a {\it
collision singularity } if there exists  $\mbox
{{\mbox{\rv}}}^{\ast}\in\cal C$ such  that $\mbox{{\mbox{\rv}}}(t)
\rightarrow \mbox{{\mbox{\rv}}}^{\ast}$, when   $t\rightarrow
c$. Otherwise, the singularity is called a
singularity without collision. Here is our first result.\\

\noindent {\bf Theorem A.} {\it All
singularities in the fixed homogeneous circle problem are collision singularities.}\\

In what follows we consider the fixed homogeneous circle   {$\,\cal C$}
contained in the $xy$-plane and centered at the origin. Also, by rescaling we can consider
the circle with radius equal to one (see section 1.3). Then the mass $M$ of $\cal C$ is given by
$M\, =\, 2\pi\lambda$.
A quick examination of the problem shows that
(see section 1.2 for more details) the {\it $z$-axis}, the {\it horizontal plane}
 (i.e. the $xy$-plane, which contains the circle)
and any {\it vertical plane} (i.e any plane that contains the $z$-axis) are invariant
subspaces of our problem. Also, any {\it radial line} (i.e.
any line in the $xy$-plane that passes through the origin)
is an invariant subspace. Here by an invariant subspace $\Lambda$ we mean that
any solution that begins tangentially in $\Lambda$ is totally contained in
$\Lambda$. The restriction of the problem to the $z$-axis is a
one dimensional problem which is not difficult to treat. In fact it is a special case
of a problem studied by Sitnikov (see \cite{BLO}, \cite{M}). Our two next results discuss the
dynamics restricted to the horizontal plane. This problem is a central force problem.
For $r\in\R-\{ 1\}$ define the function
$V(r)=V(r{\bf u})$, where ${\bf u}$ is any vector in the $xy$-plane
of length one (note that we are using the same letter $V$ for two different functions).
The dynamics restricted to the horizontal plane has two cases: inside the circle and outside the circle.
It is a classical result that using the polar coordinates
$(r,\theta)$ of a point $(x,y)$ in the horizontal plane it can be proved that
system \ref{00} in this case is equivalent to:
{\footnotesize$
\stackrel{..}{{r}}\,=\, \displaystyle\frac{K^{2}}{{r}^{3}}-\frac{d}{dr}V(r),\,\,
\stackrel{.}\theta\,=\, \displaystyle\frac{K}{{r}^2}$}
(see \cite{G}). Here the constant $K$ is the angular momentum of the solution.
Hence this problem can be reduced in a canonical way
to a problem with one degree of freedom (the variable $r(t)$) and once we know
$r(t)$, we can obtain $\theta(t)$ by integration.
Note that $K=0$ if and only if the solution lies in a radial line (i.e. $\theta (t)$ is constant).
Equation {\footnotesize$\stackrel{..}{{r}}\,=\, \displaystyle\frac{K^{2}}{{r}^{3}}-\frac{d}{dr}V(r)$}
is equivalent to

{\footnotesize\begin{equation}
\stackrel{..}{r(t)} \,=\,-\frac{d}{dr} U(r(t)), \label{01}
\end{equation}}

\noindent where  $U$ denotes the {\it effective potential}
 {\footnotesize$U({r})=\displaystyle\frac{K^{2}}{2{r}^{2}}+V({r})$}.
Note that $U(r)$ depends on $K$, that is for each $K$ we may have different functions $U(r)$
(maybe we should write $U_K(r)$ but we do not want to complicate our notation).
Hence to determine the dynamics we have to know the behavior of the effective potential
for all values of $K$.
Note that outside the circle the effective potential is defined for $1<r<\infty $ and
inside the circle the effective potential is defined for $0< r <1 $ (for $K\neq 0$)
and for $-1< r <1 $ (for $K= 0$).\\

\noindent {\bf Theorem B.} {\it  The effective potential $U(r)$ inside the
circle has the following properties:}
\begin{enumerate}
\item[{i.}]{\it   $\frac{d}{dr}V(r)<0$, for $0<r<1$. Hence $\frac{d}{dr}U(r)<0$, for $0<r<1$.}

\item[{ii.}]{\it   $lim_{r\rightarrow 1^{-}}U(r)=-\infty$.}

\item[{iii.}]{\it For $K\neq 0$ we have $\lim_{r\rightarrow 0^{+}}U(r)=\infty$.}

\item[{iv.}]{\it  For $K=0$ we have that $U(r)$, $-1<r<1 $, is even.}
\end{enumerate}
\vspace{.2in}

\noindent The following are direct consequences of the Theorem above.\\

\noindent {\bf 1.} The graphs of $U(r)$ for the cases $K=0$ and $K\neq 0$
have the following form.

\newpage

\begin{figure}[!htb]
\centering
 \includegraphics[scale=0.5]{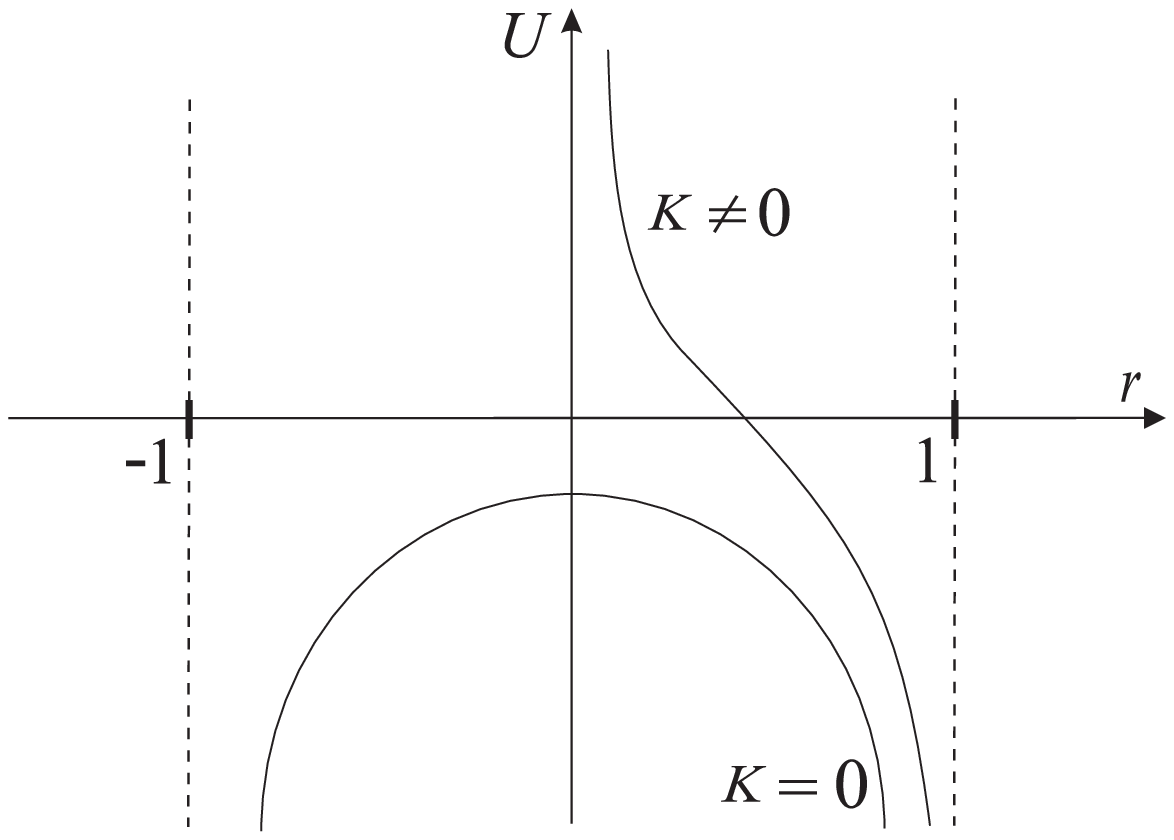}
 \caption{\scriptsize{Graph of the effective potential $U$ inside the circle.}}
 \end{figure}

\noindent Note that for $K=0$,  $U(r)$, $-1<r<1$, has an absolute
maximum value $E_0=U(0)=V(0)=-M=-2\pi\lambda$ , at $r=0$.
Hence $r\equiv 0$ or, equivalently ${\mbox{\rv}}(t)=0$, $t\in \R$,
is an equilibrium solution.\\

\noindent {\bf 2.} By (i) of the Theorem, the force inside the circle is repulsive.\\

\noindent {\bf 3.} As mentioned before, for $K=0$ a solution ${\mbox{\rv}}(t)$ (with polar coordinates
($r(t), \theta(t)$) stays in a radial line, i.e. $\theta$ is constant.
Without loss of generality we can assume that this radial line is the $x$-axis, that is $\theta=0$. Hence
${\mbox{\rv}}(t)=(x(t),0,0)$, where $x(t)=r(t)$. Let
$E=E(x,\dot x )$ denote the energy of  ${\mbox{\rv}}(t)=(x(t),0,0)$. It follows that the phase portrait
for $K=0$ (i.e the level curves of $E=E(x,\dot x )$ in the $x\dot x$-plane) has the following form:
\begin{figure}[!htb]
 \centering
 \includegraphics[scale=0.5]{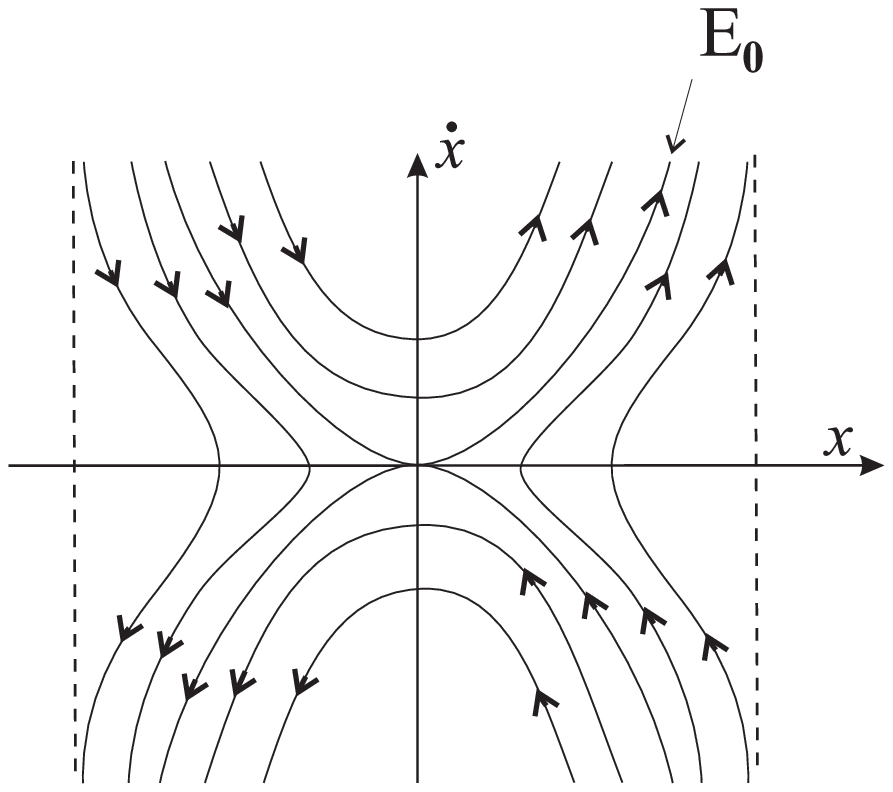}
\caption{\scriptsize{Phase portrait for $K=0$.}}
 \end{figure}

Consequently, depending on the energy the particle behaves in the following way.
If $E<E_0$ the particle comes from the point $p=(\pm 1,0,0)$ in the circle, stops before reaching
the origin and then turns back to the point $p$.
If $E>E_0$ the particle comes from the point $p$ in the circle, and goes all the way to the point $-p$.

\begin{figure}[!htb]
    \centering
\includegraphics[scale=0.5]{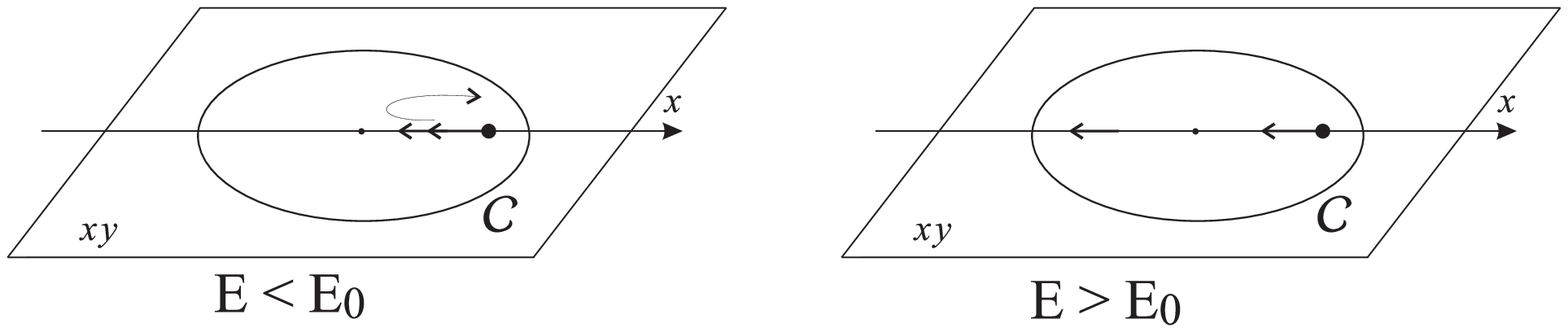}
\end{figure}

If $E=E_0$ we have three cases: the particle comes from the point $p$ in the circle
and converges to the origin (in infinite time), or comes from the origin
(in infinite time) and converges to $p$, or stays at rest in the origin.

\begin{figure}[!htb]
 \begin{minipage}[b]{0.28\linewidth}
 \includegraphics[width=\textwidth]{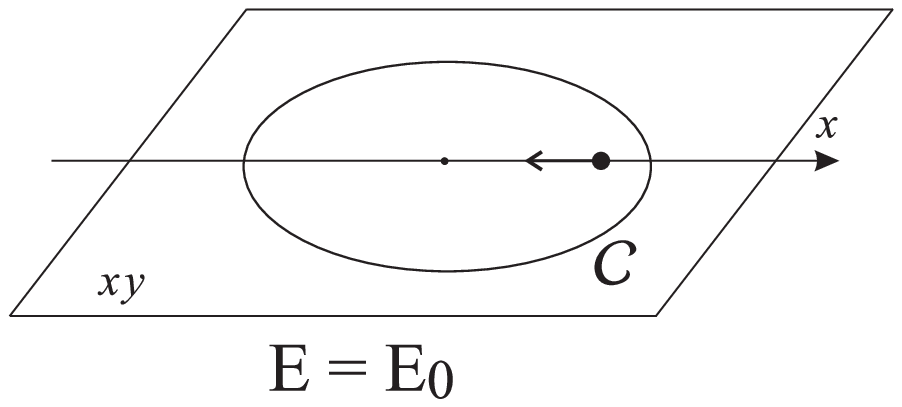}
 \end{minipage}  \hfill
 \begin{minipage}[b]{0.28\linewidth}
 \includegraphics[width=\textwidth]{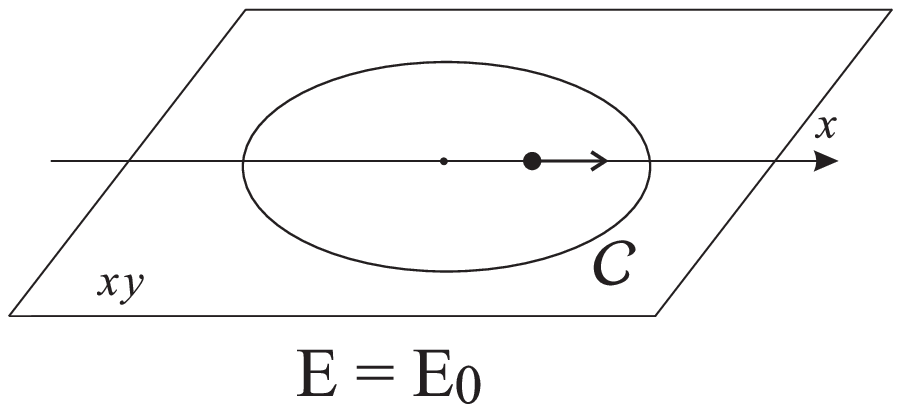}
 \end{minipage}  \hfill
 \begin{minipage}[b]{0.28\linewidth}
 \includegraphics[width=\textwidth]{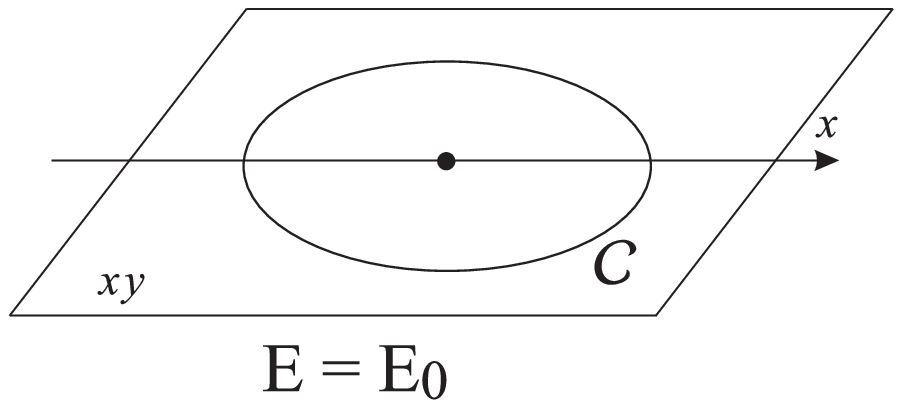}
 \end{minipage}  \hfill
 \end{figure}

\noindent {\bf 4.} For $K\neq 0$ we have that a solution ${\mbox{\rv}}(t)$
 never passes through the origin. By Theorem B, this solution comes from
the circle and returns to it (in finite time, by Theorem A). Note
that there is a unique $t_0$ such that ${\mbox{\rv}}(t_0)$ is the
closest point to the origin. Since $\dot r (t_0)=0$ we have that
$\dot {\mbox{\rv}} (t_0)$ is perpendicular to ${\mbox{\rv}}(t_0)$.
We say that ${\mbox{\rv}}(t)$ is {\it normalized } if
${\mbox{\rv}}(t_0)$ lies in the positive $y$-axis. Note that, by
symmetry, any solution (inside the circle) with $K\neq 0$ can be
obtained from a normalized one by a rotation.
(See Fig. 0.3 below.)\\

The following two Propositions complement the results given by Theorem B, for the case of a solution
inside the circle.\\

\noindent {\bf Proposition 1.} {\it Let ${\mbox{\rv}}(t)$ be a
solution of \ref{00}, contained in the horizontal plane and
inside the circle. Let $(a,b)$ be its maximal interval
of definition. If ${\mbox{\rv}}(t)$ approaches the circle as
$t\rightarrow b^-$ or $t\rightarrow a^+$, then }
\begin{enumerate}
\item[{(1)}] {\it $b< \infty $ or $ -\infty <a$, respectively. Hence, by Theorem A, the solution collides
to a point in the circle,}

\item[{(2)}] {\it when $\mbox{\rv} (t)$ converges
to the circle the speed converges to infinity and the velocity $\dot {\mbox{\rv}} (t)$
becomes orthogonal to the circle.}
\end{enumerate}
\vspace{.2in}

It follows from this Proposition that the intervals of definition $(a,b)$ of solutions inside the circle
are as follows: (1) for $K\neq 0$ we have that $a$ and $b$ are both finite
(2) for $K=0$ and $E\neq E_0$ we also have
$a$ and $b$ are finite (3) for $K=0$ and $E=E_0$ either $a=-\infty$ and $b=\infty$ or exactly one
of $a$ and $b$ is finite.\\

\noindent {\bf Proposition 2.}
{\it Let ${\mbox{\rv}}(t)$ be a solution of \ref{00}, contained in the horizontal plane,
inside the circle and with angular momentum $K\neq 0$. If ${\mbox{\rv}}(t)$ is normalized then
the trace of the curve ${\mbox{\rv}} (t)$ is the graph of
an even convex function $y=y(x)$.}

\begin{figure}[!htb]
 \centering
 \includegraphics[scale=0.7]{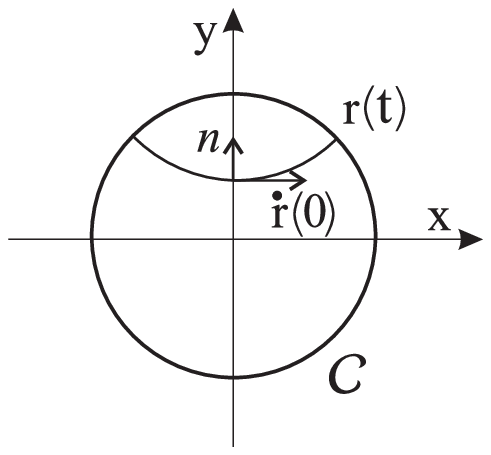}
\caption{\scriptsize{Trace of the normalized solution $\mbox{\rv}(t)$.}}
 \end{figure}

We now state the results concerning the effective potential $U(r)$ for the case
$r>1$, that is outside the circle.\\

\noindent {\bf Theorem C.} {\it  The effective potential $U(r)$, $r>1$ (outside the
circle) has the following properties:}
\begin{enumerate}
\item[{i.}]{\it $\frac{d}{dr}\, V(r)>0$.  Hence $\frac{d}{dr}(U(r)-\frac{K^2}{2r^2})>0$.}

\item[{ii.}]{\it   $lim_{r\rightarrow 1^{+}}\, U(r)=-\infty$.}

\item[{iii.}]{\it $lim_{r\rightarrow \infty}\, U(r)=0.$}
\end{enumerate}

\noindent {\it Also there are $K_0>0$, $r_0$, $1<r_0<2$,
and continuous functions $r_1=r_1(|K|)$, $r_2=r_2(|K|)$ defined for $|K|\geq K_0$,
 with $1<r_1\leq r_0 \leq r_2$ such that:}
\begin{enumerate}
\item[{iv.}] {\it For $K\in (-K_0,K_0)$, $U(r)$ does not have critical points.
Hence from ii. $\frac{d}{dr}\, U(r)>0$ and $U(r)$ is an increasing
function.}

\item[{v.}]{\it  For $K\notin (-K_0,K_0)$ the critical
points of $U(r)$ are exactly $r_1=r_1(|K|)$, $r_2=r_2(|K|)$.}

\item[{vi.}] {\it $\frac{d}{dr}\, U(r)>0$ for $r<r_1$ and $r>r_2$. Also $\frac{d}{dr}\, U(r)<0$
for $r_1<r<r_2$.}

\item[{vii.}] {\it $r_1(K_0)=r_2(K_0)=r_0$.}

\item[{viii.}] {\it $r_1$ is a decreasing function and $lim_{K\rightarrow \infty} r_1=1$}.

\item[{ix.}] {\it $r_2$ is an increasing function and $lim_{K\rightarrow \infty } r_2=\infty$}.

\item[{x.}] {\it $lim_{K\rightarrow \infty } U(r_1(K))=\infty$ and
$\,\,lim_{K\rightarrow \infty } U(r_2(K))=0$}.
\end{enumerate}
\vspace{.2in}

\noindent {\bf Remark.} The special value $r_0$ does not depend on the density $\lambda$ (or on the mass $M$)
of the circle (see the proof of Theorem C).
 Since $K_0^2=r_0^3\frac{d}{dr}V(r_0)$, the value $K_0$ does depend on $\lambda$.
It is interesting to compare this with Corollary \ref{1.3.2}.\\

\noindent It follows from Theorem C that the graphs of $U(r)$ have the following form:

\begin{figure}[!htb]
\centering
 \includegraphics[scale=0.5]{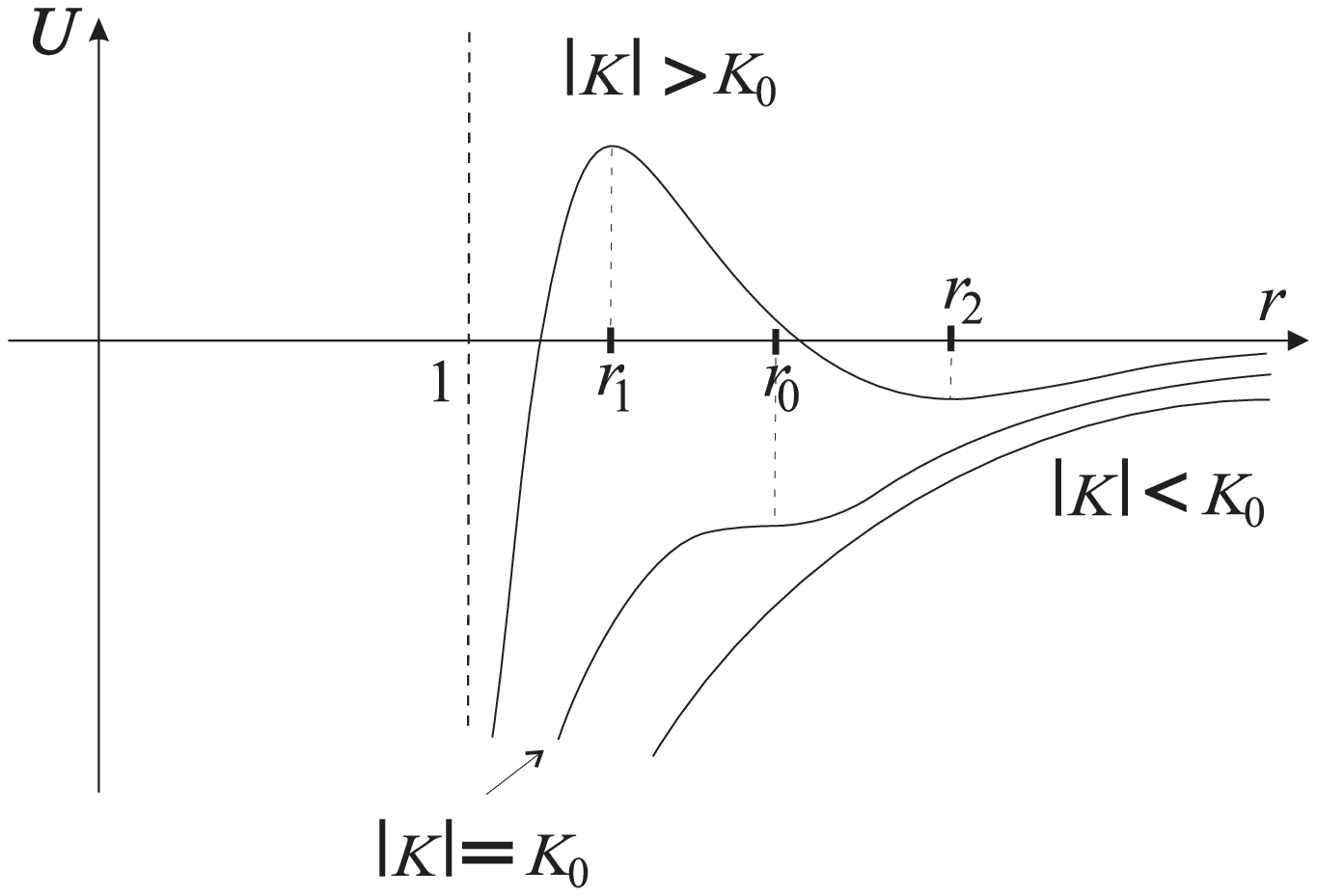}
 \caption{\scriptsize{Graph of the effective potential $U$ outside the circle.}}
 \end{figure}

\noindent Note that for $K\notin [-K_0,K_0]$ the function $U(r)$
has a local maximum at $r_1$ and a local minimum at $r_2$. By ({\it x.}) of Theorem C,
the local maximum value $U(r_1)$ at $r_1$ is eventually
(as $K\rightarrow \infty$ and $r_1 \rightarrow 1^-$) a global maximum and $U(r_1)\rightarrow +\infty$.
The following are direct consequences of Theorem C.\\

\noindent {\bf 1.} From (i) of the Theorem, the force outside the circle is attractive.\\

\noindent {\bf 2.} We discuss the stability of circular solutions.
Here by stability we mean the following.
 A circular solution ${\mbox\rv}(t)$ in the horizontal plane is {\it stable} if
 its reduced solution $r(t)=\|{\mbox\rv}(t)\|=$ constant is a {\it stable} equilibrium position
 of the reduced problem $\stackrel{..}{r(t)} \,=\,-\frac{d}{dr} U(r(t))$.
It is easy to see that a
 circular solution ${\mbox\rv}(t)$ is stable iff the trace of solutions with initial
conditions close to the initial conditions of ${\mbox\rv}(t)$ stay
close to the trace of ${\mbox\rv}(t)$. Recall that for a fixed
$K$, a circular solution with radius $r$ has momentum $K$ if and
only if $r$ is a critical point of $U(r)$. Hence for $K\in
(-K_0,K_0)$ there are no circular solutions with momentum $K$. For
every $K\notin [-K_0,K_0]$ there are exactly two circular
solutions with radii $r_1(|K|)$ and $r_2(|K|)$. The circular
solution with radius $r_1$ is not stable and there are ``spiral''
solutions approaching it from the inside and the outside. The
circular solution with radius $r_2$ is stable and has no
``spiral'' solutions approaching it.
 For $K=K_0$ there is exactly one circular solution and it has radius $r_0$.
This circular solution is unique among all circular solutions. It is not stable and
have spiral solutions approaching it only from the inside. Also it can be approximated by stable
circular orbits from the outside (with varying momenta $K$). It follows that a circular
solution with radius $r$ is stable if and only is $r>r_0$.
Note that a non-circular solution with momentum $K_0$ either collides to the circle,
escapes to infinity
or is the solution that approaches the circle (in infinite time).\\

\noindent {\bf 3.} Let $E=E(r,\dot r )$ denote the energy of  $r(t)=\|{\mbox\rv}(t)\|$.
 Using Theorem C we can deduce the phase
 portraits. We have three cases: $K\in(-K_0,K_0)$, $K=K_0$, $K\notin [-K_0,K_0]$.
 For the first two cases the phase portraits are shown in figures 0.5 and 0.6 respectively,
 and the solutions have the following behavior:
 \vspace{0.2cm}

 \noindent(i) $K\in(-K_0,K_0)$. In this case all solutions either
 collide to the circle or escape to infinity.
 \vspace{0.3cm}

 \noindent(ii) $K=K_0$. In this case all solutions (different from the unique
circular solution of radius $r_0$) either collide to the circle or
escape to infinity or converge to the circular solution with radius
$r_0$ (in infinite time).

\begin{figure}[!htb]
\hspace{2cm}
 \begin{minipage}[b]{0.3\linewidth}
 \includegraphics[width=\textwidth]{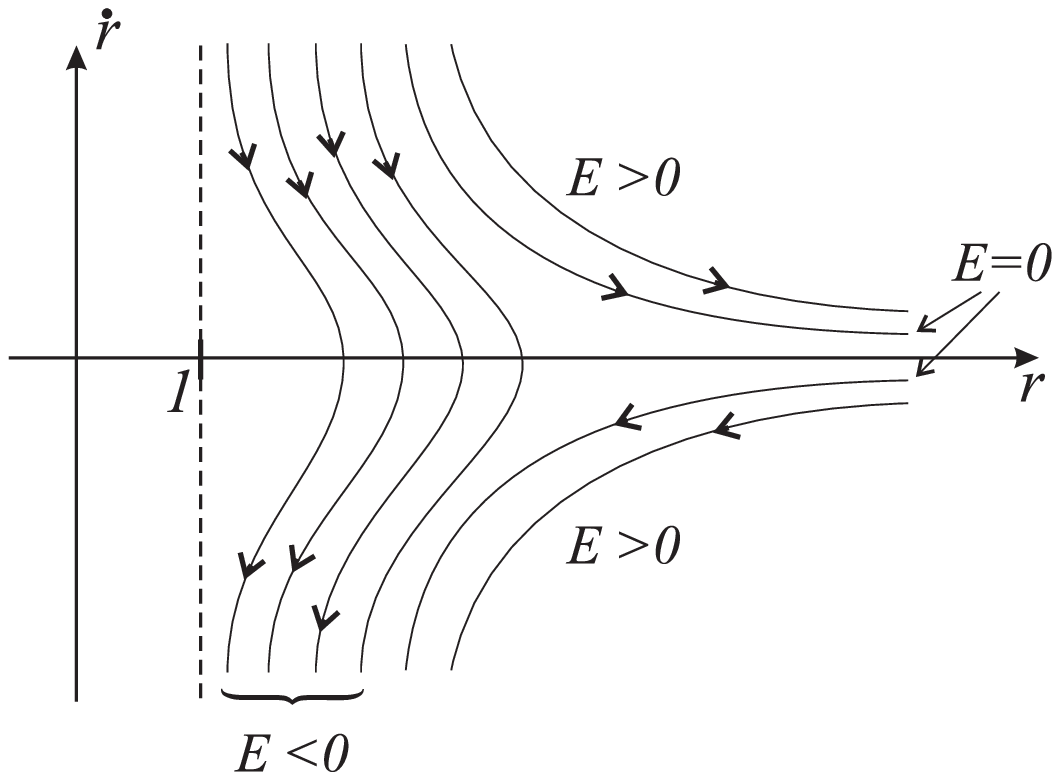}
 \caption{\scriptsize{$K\in (-K_0,K_0)$.}}
 \end{minipage}  \hfill
 \begin{minipage}[b]{0.3\linewidth}
 \includegraphics[width=\textwidth]{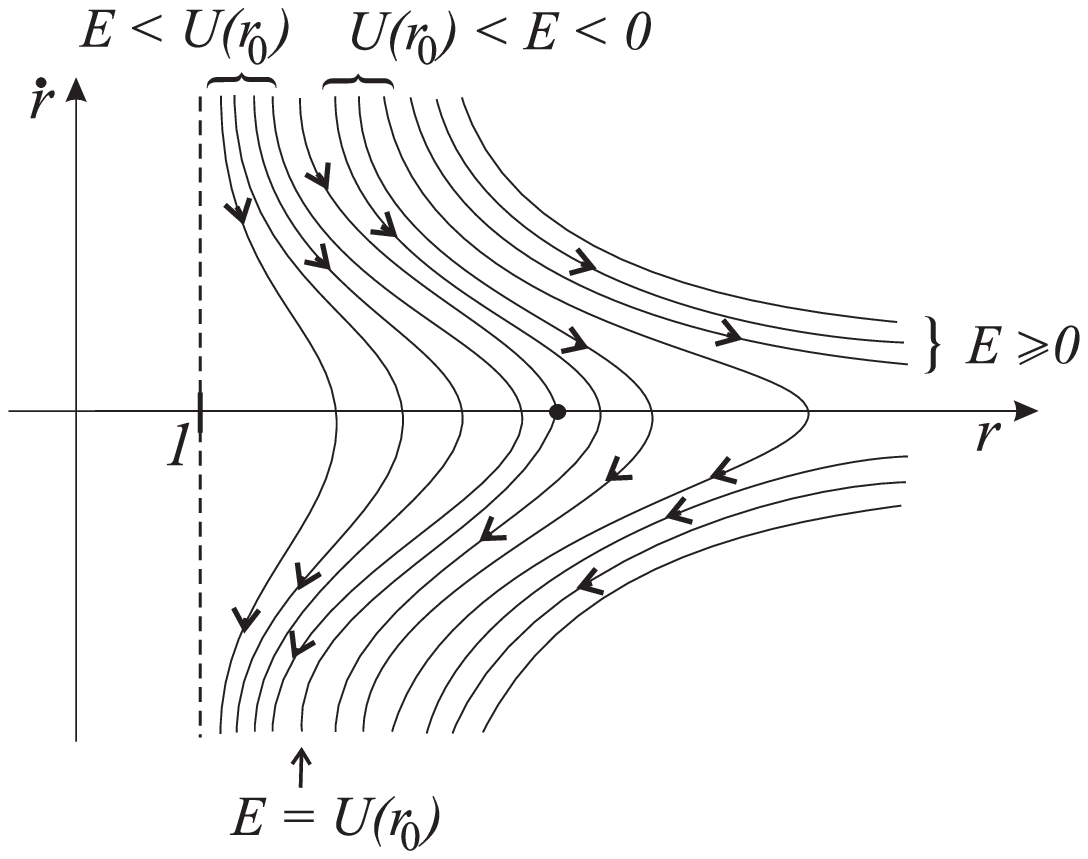}
 \caption{\scriptsize{$K=K_0$}}
 \end{minipage}
  \hspace{2cm}
 \end{figure}

 \noindent (iii) We analise now the case $K\notin [-K_0,K_0]$.
Let $\overline{E}=U(r_1)$. We have also three cases:
 $\overline{E}>0, \overline{E}<0, \overline{E}=0$ (see figures 0.7, 0.8, 0.9).
Note that, by ({\it vi.}), ({\it vii.}) and ({\it x.}) of Theorem C
all three cases do happen. For $\overline{E}>0$ we have (see figure 0.7):
 In  region I the solutions are bounded, and converge to
the circle. In  region II the solutions are unbounded,  come from
infinity and converge   to the circle. In  region III the solutions
are unbounded, come from  the circle and escape to infinity. In
region IV the solutions are unbounded,  approach the circular
solution (with radius $r_1$) and escape to infinity. In  region V
the solutions are bounded, and stay close to the circular solution
(with radius $r_2$).

 In figure 0.7 the solution represented by $\gamma _{1}$ comes from the circle and converges
to the unstable circular solution (with radius $r_1$).
 The solution represented by $\gamma _{2}$ comes from the circular solution
(with radius $r_1$) and converges to the circle.
 The solution represented by $\gamma _{3}$ comes from the circular solution
(with radius $r_1$) and escapes to infinity.
The solution represented by $\gamma _{4}$ comes from infinity
and converges to the circular solution (with radius $r_1$).\vspace{0.2cm}

For $\overline{E}=0$ we have also three cases: $E<\overline{E}$, $E=\overline{E}$, $E>\overline{E}$.
 We have the following analysis  (see figure 0.8 below):
For $E< \overline{E}$
in  region I the solutions are bounded, and converge to the circle;
in  region IV the solutions are bounded and stay close to the circular solution (with radius $r_2$).
For $E>\overline{E}$ we have:
In  region II the solutions are unbounded, come from  the
circle and escape to infinity.
In region III the solutions are unbounded,  come from infinity and converge to the circle.
For $E=\overline{E}$ we have:
The solution represented by $\gamma _{1}$ comes from the circle and converges
to the unstable circular solution (with radius $r_1$).
 The solution represented by $\gamma _{2}$ comes from the circular solution
(with radius $r_1$) and converges to the circle.
 The solution represented by $\gamma _{3}$ comes from the circular solution
(with radius $r_1$) and escapes to infinity.
The solution represented by $\gamma _{4}$ comes from infinity
and converges to the circular solution (with radius $r_1$).\vspace{0.2cm}

For $\overline{E}<0$ we have four cases: $E<\overline{E}$, $E\geq
0$, $\overline{E}<E<0$, $E=\overline{E}$ (see figure 0.9 below).
 For the case $E< \overline{E}$ in  region
I  the solutions are bounded and converge
to the circle; in  region IV the solutions are bounded and stay
close to the circular solution (with radius $r_2$). For $E\geq 0$
we have: In  region II the solutions are unbounded, come from  the
circle and escape to infinity. In region III the solutions are
unbounded,  come from infinity and converge to the circle. For
$\overline{E}<E<0$ the solutions are bounded, come from  the
circle, approaching the circular solutions and return converging
to the circle. Finally, for $E=\overline{E}$ we have: The solution
represented by $\gamma _{1}$ comes from the circle and converges
to the unstable circular solution (with radius $r_1$).
 The solution represented by $\gamma _{2}$ comes from the circular solution
(with radius $r_1$) and converges to the circle.
 The solution represented by $\gamma _{3}$ comes from the circular solution
with radius $r_1$, approach the circular solution with radius
$r_2$, and return converging to the circular solution
with radius $r_1$.

\begin{figure}[!htb]
 \begin{minipage}[b]{0.32\linewidth}
 \includegraphics[width=\textwidth]{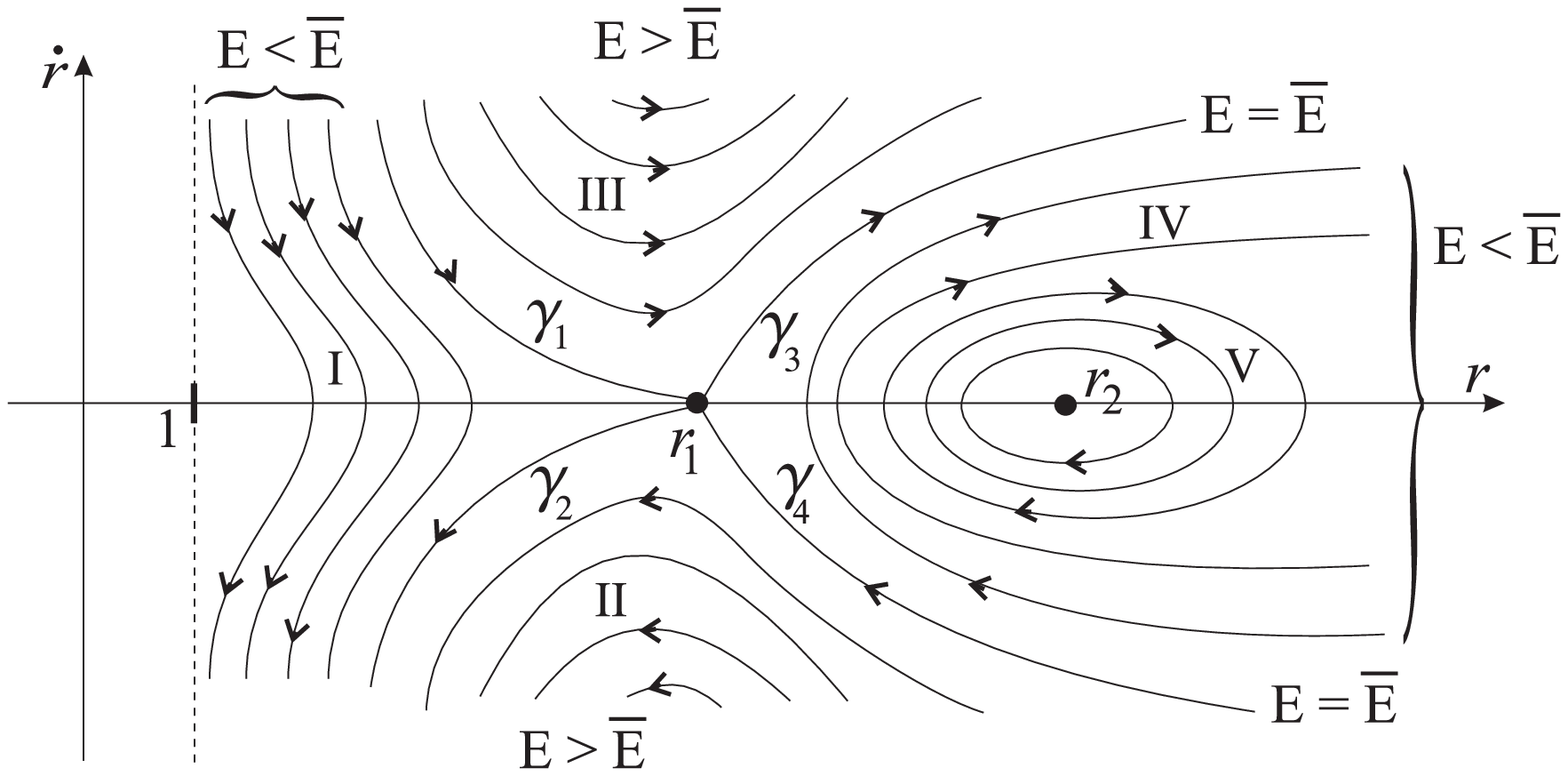}
 \caption{\scriptsize{$\overline{E}=U(r_1)>0.$}}
 \end{minipage}  \hfill
 \begin{minipage}[b]{0.32\linewidth}
 \includegraphics[width=\textwidth]{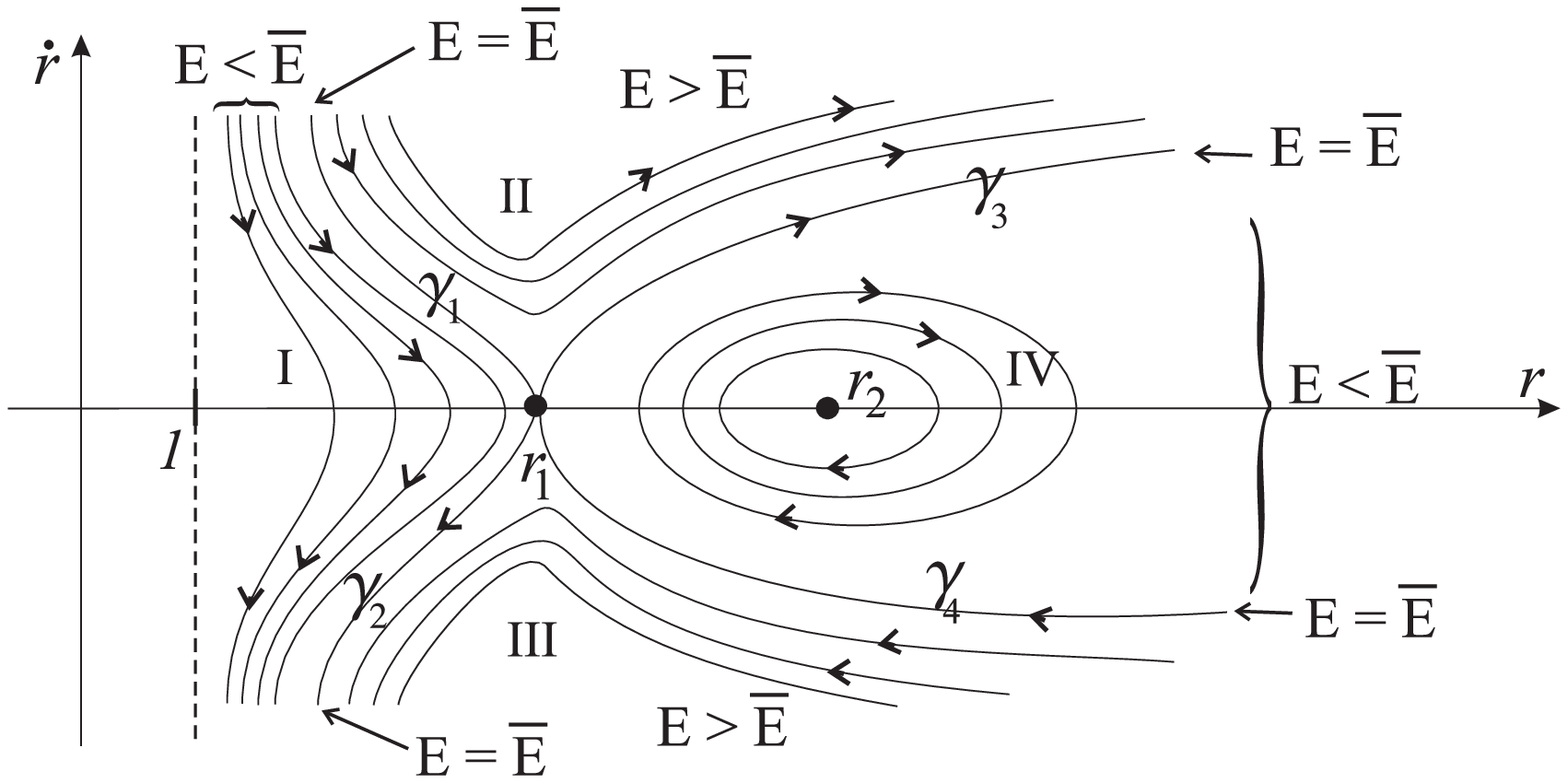}
 \caption{\scriptsize{$\overline{E}=U(r_1)=0.$}}
 \end{minipage}  \hfill
 \begin{minipage}[b]{0.3\linewidth}
 \includegraphics[width=\textwidth]{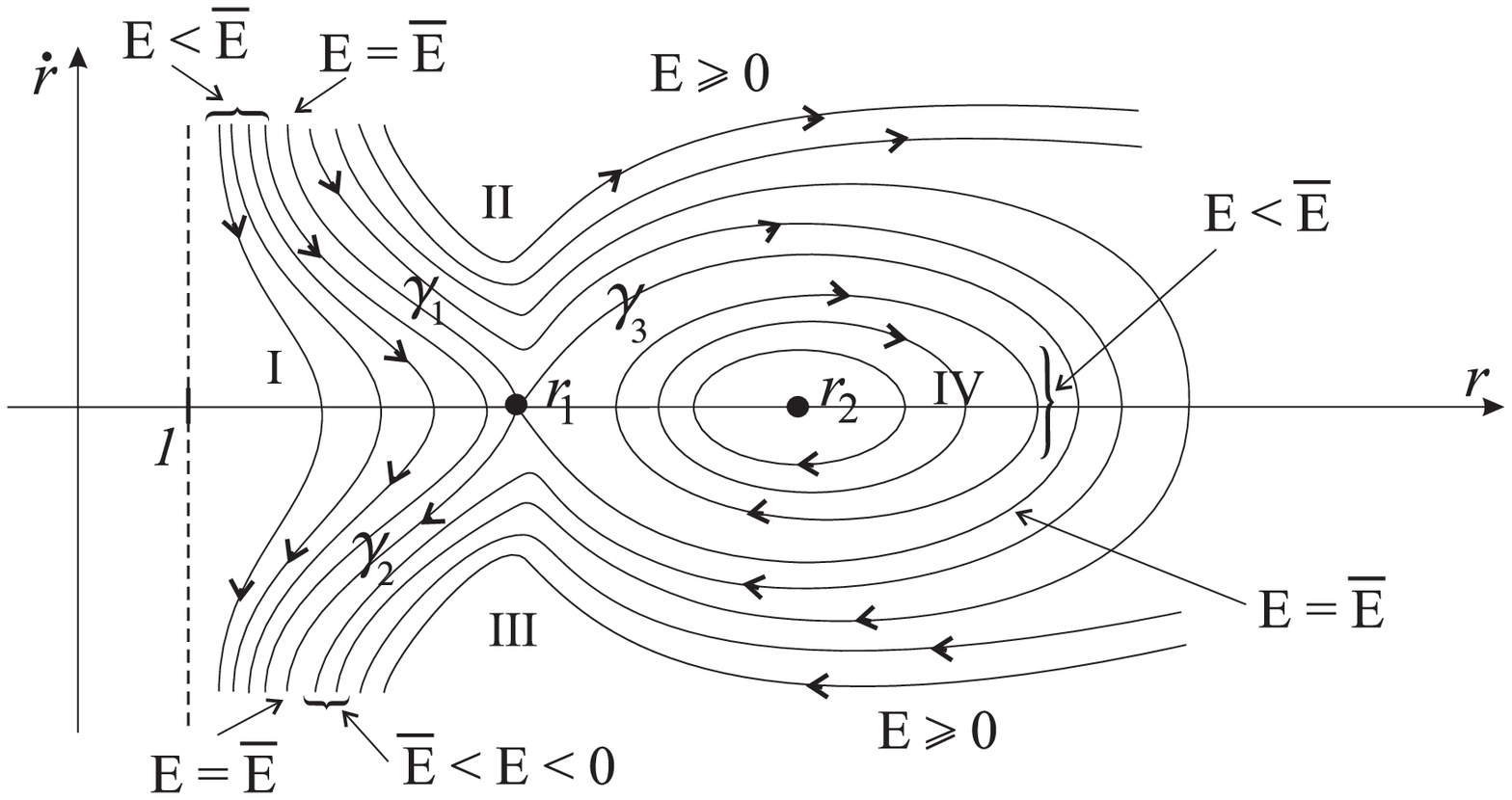}
\caption{\scriptsize{$\overline{E}=U(r_1)<0.$}}
 \end{minipage}  \hfill
 \end{figure}

\vspace{0.4cm}

We mentioned above that a circular solution of radius $r$ in the
horizontal plane is stable (among solutions in the horizontal
plane) if and only if $ r>r_0$. The following question arises: are
these solutions stable among all solutions in $\R^3$? The answer is affirmative, but
before we state this result we have some comments. Since our
original problem (in $ \R^3$) is invariant by rotations around the
$z$-axis, the problem can be reduced in a canonical way to a
problem with two degrees of freedom (see section 5). In fact,
every solution ${\mbox \rv}(t)=(x(t),y(t),z(t))$  can be written
as $\mbox{\rv}(t)=(r(t)cos\varphi (t), r(t)sin\varphi(t),z(t))$,
i.e the cylindrical coordinates of ${\mbox \rv}(t)$ are
$(r(t),\varphi(t),z(t))$, and $(r(t),\varphi(t),z(t))$ satisfy a system of
equations (system (\ref{R4})). Moreover, once we know $r(t)$ we can obtain
$\varphi(t)$ by integration. We say that $(r(t),z(t))$ is the {\it
canonical projection solution} of ${\mbox \rv}(t)$. That is, the
canonical projection solution is obtained from the solution by the
map $(x,y,z)\mapsto (r,z)$, $r=\sqrt{x^2+y^2}$. We have that
$(r(t),z(t))$ satisfy the first two equations of (\ref{R4}). Note that under
this projection the circle $\cal C$ projects to a point $x_{\cal
C}$ with $rz$ coordinates (1,0).
We extend the definition of a circular solution (for central force problems)
to this context in the following way:
 We say that a solution is {\it
circular} if its canonical projection solution $(r(t),z(t))$ is
an equilibrium position of the system formed by the first two equations of (\ref{R4}).
Also, we say that a circular solution is {\it stable} if its canonical projection solution
is a {\it stable} equilibrium position of the system formed by the first two equations of
(\ref{R4}). It is easy to see that a circular solution ${\mbox\rv}(t)$ is stable iff the trace
of solutions with initial
conditions close to the initial conditions of ${\mbox\rv}(t)$ stay
close to the trace of ${\mbox\rv}(t)$.\\

\noindent {\bf Proposition 3.} {\it
 All circular solutions lie in the horizontal plane. Moreover, a circular solution
(in the horizontal plane) of radius $r$ is stable in $\R^3$ if and only if $r>r_0$,
where $r_0$ is as in Theorem C.}\\

We present two more results. Fix $\lambda >0$. For
$\epsilon >0$ denote by $W({\mbox{\rv}},\epsilon)$ the potential
function induced by a the fixed homogeneous circle in the
$xy$-plane centered at $(\frac{1}{\epsilon},0,0)$ with radius
$\frac{1}{\epsilon}$ and density $\lambda$. We write $W(x,z;\epsilon)$ for the
restriction of this potential to the $xz$-plane. Define $\nabla
W(x,z;0):=64\,\lambda\frac{(x,z)}{x^{2}+z^{2}},$ $(x,z)\neq
(0,0).$ Let $\nabla W(x,z;\epsilon)$ be the gradient of
$W(x,z;\epsilon)$ with respect to the variables $x$ and $z$.
 Also, let $A=\{\,(x,z;\epsilon)\,;\,
(0,0)\neq \,(x,y)\neq (\frac{2}{\epsilon},0)\,\}$.
Thus $\nabla W(x,z;\epsilon)$ is defined on $A$. \\

\noindent {\bf Proposition 4.} {\it  $\nabla W$ is continuous on   $A$.}\\

Note that $\nabla W(x,z;0)=64{\lambda}\nabla ln
(\sqrt{x^{2}+z^{2}})$ and recall that $2{\lambda}\,ln
(\sqrt{x^{2}+z^{2}})={\lambda}ln (x^{2}+z^{2})$ is the  potential
induced by the infinite wire (with constant density $\lambda$ and
infinite mass) orthogonal to $xz$-plane intersecting the
$xz$-plane at the   origin. Hence by Proposition 4 the problem of
the fixed homogeneous circle with large radius, constant density
$\lambda$, and centered at $(\frac{1}{\epsilon},0,0)$) can be
regarded as a perturbation of the problem of the infinite
homogeneous straight wire with density $32\lambda$. (To obtain the
same density $\lambda$, instead of $32\lambda$, it is enough to
take circles with radius $\frac{1}{2\epsilon}$.) This result is a
key element in the proof of the existence of periodic
orbits near the circle (see \cite{AO}).\\

There are a few textbooks where it is claimed that the potential
function of fixed homogeneous circles with large radii ``converge''
to the potential function of the infinite straight wire, but we have
not found satisfactory proofs of these claims. In fact, it seems
that in order to avoid the problems with the infinities, we have to
work with the gradients, instead of the potentials. This necessity
seems to be related to the fact that if we try to calculate the
potential of the infinite straight wire by integrating
$\int^{+\infty}_{-\infty} \frac{dx}{\| p-(x,0,0)\|}$ we obtain
infinity, but the integral for the gradient does converge (or can be
computed easily from Gauss formula). \\

Finally we present our last result. The potential $V$ of our
problem at $P$ can be written as $V(P) =-\frac{M}{\sigma (D,d)}$,
where $\sigma (D,d)$ is the arithmetic-geometric mean of $D>0$ and
$d>0$  (see Remark 1.1). Here $d, D$ are the maximum and minimum
distances $D=D(P)$, $d=d(P)$ of the point $P$ to the circle.
 As we mentioned above, this formula of $V(P)$  was given by Gauss.
In our last result we give a formula for $(\frac{\partial}{\partial
D}V,\frac{\partial}{\partial d}V)$
  in terms of $D, d$ and their successive geometric and arithmetic means $d_j, D_j$. \\

\noindent{\bf Proposition 5.} {\it We have the following formulas:
$$\frac{\partial}{\partial D}V(D,d)=\frac{\chi -1}{D} V(D,d),\,\,\,\,\,\,\,
\frac{\partial}{\partial d} V(D,d)= \,-\frac{\chi}{d}V(D,d),$$

\noindent where $\chi=\sum_{n=1}^{+\infty} \frac{1}{2^{n}}
      \frac{ d_{n-1}}{D_{n}}\prod_{j=1}^{n-1}\left(
      1-\frac{d_{j-1}}{D_{j}}\right)$.}\\

\vspace{0.4cm}

This paper has 7 sections. In section 1 we give some preliminary
results that will be needed later. These include some basic facts
about symmetry, elementary properties of the potential function as
well as some identities related to the potential function. In
sections 2, 3 and 5 we prove Theorems A, B and C, respectively. At
the end of section 5 we prove Proposition 3 and also prove
that the only equilibrium position is the origin. By symmetry it
is clear that the origin is an equilibrium position but it is not
too trivial to prove that it is the only one. In section 4 we
prove Propositions 1 and 2. Finally in sections 6 and 7 we prove
Propositions 4 and 5, respectively.\\

 \section{ Preliminaries.}
\subsection{The Potential.}

\begin{wrapfigure}[8]{o}{4.5cm}
 \centering
 \includegraphics[width=4cm]{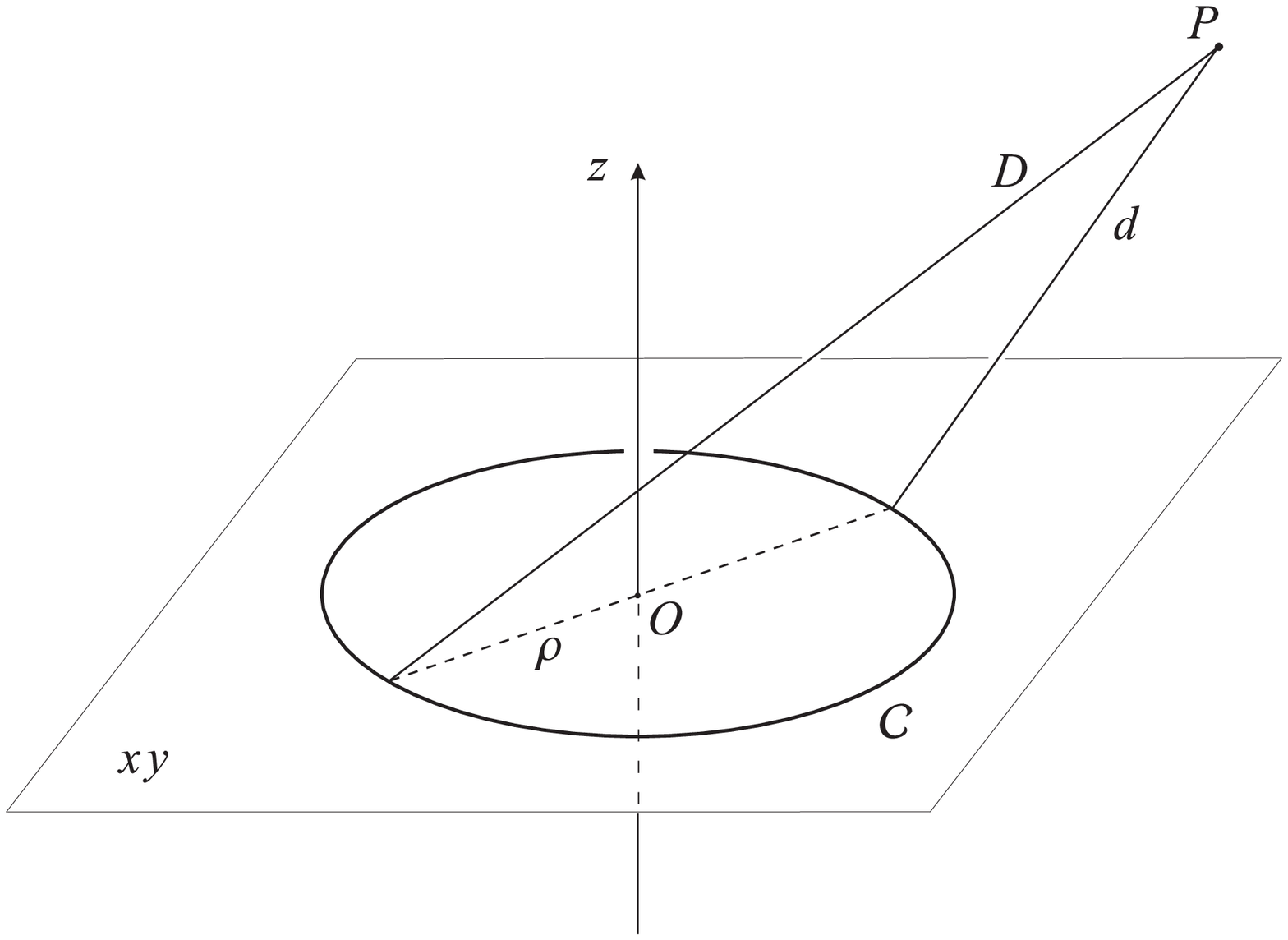}
 \caption{\scriptsize{Fixed circle centered at the origin. }}
 \end{wrapfigure}

As before consider the fixed homogeneous circle   {$\,\cal C$} with constant
density $\lambda$ and radius $\rho$,
contained in the $xy$-plane, and centered at the origin.
The mass $M$ of $\cal C$
is $2\pi \lambda\rho$.
We denote by
$P=(x,y,z)$ the coordinates of the position of a  particle in
 $\R^{3}\setminus\cal C$ and let $V(P)$ be the potential at $P$ induced by $\cal C$.
Let $D=D(P)$  and $d=d(P)$ denote the maximum and minimum distances
from the particle $P$ to the circle. We have the following
expressions for these distances: $D^{2} = (\sqrt{x^2 +y^2}+\rho)^2
+z^2$ and $d^{2} = (\sqrt{x^2 +y^2}-\rho)^2 +z^2.$ Also we have that
(see \cite{Mc}, p.196)
\begin{equation}
V(P) = -4\lambda\rho \int_{0}^{\frac{\pi}{2}}\frac{d\psi}{\sqrt{d^2 cos^2 \psi
+ D^2 sin^2 \psi}}.
\label{1.3}
\end{equation}

\vspace{0.2cm}

Define $T:\{(x,y); \,x>0, y>0\}\longrightarrow\R$  by
{\small $T(d,D)= \int_{0}^{\frac{\pi}{2}}\frac{d\psi}{\sqrt{d^2 cos^2 \psi
+ D^2 sin^2 \psi}}$}.
 Then we have $V(P)=-4\lambda\rho \,T(d,D)$.
By the change of variable $\psi =
\frac{\pi}{2}-\theta$, it can be proved that
 $\,\,T$ is symmetric with respect to $d$ and $D$, that is, $T(d,D)=T(D,d)$.
 It was proved by Gauss (see \cite{Mc}, p.197-199) that $V(P)=-\frac{M}{\sigma(d,D)}$ where
 $\sigma(d,D)$ is the arithmetic-geometric mean of $d$ and $D$. Hence $T(d,D)=\frac{\pi}{2\sigma(d,D)}$.

\vspace{0.15in}

{\obs {\rm {Recall that the arithmetic-geometric mean $\sigma(m,n)$ of two positive
 numbers $m, n$ is defined in the following way.
Set $m_1 =\frac{m+n}{2}$ and $n_1 =\sqrt{mn}$  and define
inductively the sequences $\{ m_i\}$, $\{ n_i \}$ by
$m_{i+1} =\frac{m_i +n_i}{2}, \,  n_{i+1} =\sqrt{m_{i}n_{i}}.$
 We have $m>m_1 >m_2>...>n_2 >n_1 >n.$ It can be shown that $\lim m_i =\lim n_i$.
This common limit is called the {\it arithmetic-geometric mean }
$\sigma(m,n)$ of $m$ and $n$.
 Note that if $m=n,$ then $\,\sigma(m,m)=m.$ Note also that $\sigma(m,n)=\sigma(m_i ,n_i).$
 Since $T(m,n)=\frac{\pi}{2\sigma(m,n)}$, we also have $T(m,n)=T(m_i,n_i)$.}}
\label{1.1}}

\vspace{0,3cm}

 Define now $f(t)= \int_{0}^{\frac{\pi}{2}}\frac{d\theta}{\sqrt{
cos^2 \theta + t sin^2 \theta}}, \,\,\,\,\,\,0<t\leq 1.$
Then {\small{$T(d,D) = \displaystyle\frac{f \left( \frac{d^2}{D^2}\right)} {D}=
 \displaystyle\frac{\pi}{2\sigma (d,D)},$}}
$\,\,0<d\leq D$.
Hence
{\small{$V(P)=-\frac{4\lambda
\rho}{D}f\left( \frac{d^2}{D^2}\right).$}}
It follows that {\small{$ f\left(\frac{d^2}{D^2}\right)= \frac{\pi
D}{2\sigma (d,D)}$}}, $\,\,\,0<d\leq D.$ Setting $D=1$ we have  that
$f(d^2) =\frac{\pi}{2\sigma (d,1)}, \,\,\,d\leq 1.$
\vspace{0,5cm}

\noindent {\bf Remarks.}
 \vspace{0,1cm}

 \noindent (1) $f(t)= {\displaystyle \int_{0}^{\frac{\pi}{2}}}
\frac{d\theta}{\sqrt{ cos^2 \theta + t sin^2 \theta}}=
{\displaystyle \int_{0}^{\frac{\pi}{2}}} \frac{d\theta}{\sqrt{
1-(1-t) sin^2 \theta}} = K(\sqrt{1-t}\,)$,
 where $K$ is the elliptic integral of the first kind.
\vspace{0,2cm}

\noindent (2) $f(t)= \frac{\pi}{2}\left( 1 + {\displaystyle \sum_{n=1}^{+\infty}}
\left( \frac
{1\cdot 3 \cdots {(2n-1)}}{2\cdot 4 \cdots 2n}\right) ^2 (1-t)^n\right)$.
 This series can be obtained from the power series
of $K$ (see \cite{Mc}, p.197). \vspace{0,3cm}

{\lem

 (i) $f$ is  a decreasing function and   $f'$ is an increasing function.

\hspace{1,3cm} (ii) $\sigma (d,1)\geq \sqrt{d}$

\hspace{1,3cm} (iii) $f(t)\leq
\frac{\pi}{2}\frac{1}{\sqrt[4]{t}}$

\hspace{1,3cm} (iv) $|f'(t)|\leq \frac{1}{2}\frac{f(t)}{t}$

\hspace{1,3cm} (v) $|f'(t)|\leq \frac{\pi}{4}\frac{1}{t^{5/4}}$

\hspace{1,3cm} (vi) $\lim_{t\rightarrow 0^+} f(t)=+ \infty$

\hspace{1,3cm} (vii) $\lim_{t\rightarrow 0^+}t f'(t)=- \frac{1}{2}$
\label{1.6.2}}

\vspace{0,3cm}

\noindent  {\bf  Proof.}
 $(i)$
{\small$f'(t)=-\frac{1}{2} \int_{0}^{\frac{\pi}{2}}\frac{sin^{2}\theta d\theta}{{
(cos^2 \theta + t sin^2 \theta)}^{3/2}}<0,$} which implies
that $f$ is a  decreasing function.
{\small $f''(t)=\frac{3}{4} \int_{0}^{\frac{\pi}{2}}\frac{sin^{4}\theta d\theta}{{
(cos^2 \theta + t sin^2 \theta)}^{5/2}}>0$}, which  implies
that $f'$ is an  increasing function.
 Note that  $f^{(n)}(t)<0$, if  $n$ is odd, and
$f^{(n)}(t)>0$, if  $n$ is even.
$(ii)$
By the definition of  $\sigma$,  we have that  $\sigma (d,1)
\geq (\mbox{geometric mean of $d$ and}\,\, 1)=\sqrt{d}.$
$(iii)$ Setting $D=1,\,\,d^{2}=t$, and
applying $(ii)$, we obtain $(iii)$.
$(iv)$ {\small $|f'(t)|= \frac{1}{2} \int_{0}^{\frac{\pi}{2}}
\frac{sin^{2}\theta d\theta}{{(cos^2 \theta + t sin^2 \theta)}^{3/2}}
 = \frac{1}{2} \int_{0}^{\frac{\pi}{2}}\frac{1}{t}\frac{1}{{
(cos^2 \theta + t sin^2 \theta)}^{1/2}}\frac{tsin^{2}\theta }{{
(cos^2 \theta + t sin^2 \theta)} }d\theta
\,\,\leq \,\, \frac{1}{2}\frac{f(t)}{t},$}
(because $\frac{tsin^{2}\theta }{{
(cos^2 \theta + t sin^2 \theta)} } <1,$ $\,\,t>0$.)
The proof of $(v)$ follows directly from  $(iii)$ and  $(iv)$.
$(vi)$
Since $\sqrt{cos^2\theta +t sin^2\theta}\leq
\sqrt{cos^2\theta +t}\leq cos\theta +\sqrt{t}$ we have
$f(t)\geq \int_0^\frac{\pi}{2}\frac{d\,\theta}{cos\theta +\sqrt{t}}$.
A change of variable shows $\int_0^\frac{\pi}{2}\frac{d\,\theta}{cos\theta +\sqrt{t}}=
\int_0^\frac{\pi}{2}\frac{d\,\theta}{sin\theta +\sqrt{t}}.$
Hence, since $sin\theta\leq\theta$ we have
$f(t)\geq \int_0^\frac{\pi}{2}\frac{d\,\theta}{sin\theta +\sqrt{t}}\geq
\int_0^\frac{\pi}{2}\frac{d\,\theta}{\theta +\sqrt{t}}=
ln(\frac{\pi}{2}+\sqrt{t})-ln(\sqrt{t})\rightarrow\infty$ as $t\rightarrow 0^+$.
$(vii)$
 We have the following power series for
the elliptic integral of the first kind $K(k)$ (see \cite{Mc},
 p.203), $K(k) = \int_{0}^{\frac{\pi}{2}}\frac{d\theta}{\sqrt{
1-k^{2} sin^2 \theta}} =
\left[ 1+ \frac{k_{1}^{2}}{4} +\frac{9}{64}k_{1}^{4} +\ldots\right]
ln\frac{4}{k_{1}}-\left[  \frac{k_{1}^{2}}{4}
+\frac{21}{128}k_{1}^{4}+\ldots\right], $
where $k_{1}^{2} =1-k^{2}.$
Since $f(t) =K(\sqrt{1-t})$, taking $k_{1}^{2}=t$ we have
$f(t)=\left[
1+\frac{t}{4}+\frac{9}{64}t^{2}+\ldots\right]
ln\frac{4}{\sqrt{t}}
-\left[ \frac{t}{4}+\frac{21}{128}t^{2}+\ldots\right]$. It follows that
$f'(t)=\left[ \frac{1}{4}+\frac{18}{64}t+\ldots\right] ln \frac{4}{\sqrt{t}}-
\left[ \frac{1}{4}+\frac{42}{128}t+\ldots\right] +\left[ 1+\frac{t}{4}
+\frac{9}{64}t^{2}+\ldots\right] \left( -\frac{1}{2}\frac{1}{t}\right)$
therefore
$lim_{t\rightarrow 0}f'(t)t=-\frac{1}{2}$. \CaixaPreta

 \vspace{0,5cm}

\subsection{Symmetries of the potential.}

 In this section we study the symmetries of the  fixed
homogeneous circle problem and determine the invariant subspaces.
 As before, we are  considering the
fixed homogeneous circle $\cal C$, contained in the $xy$-plane
and centered at the  origin. We call the
$xy$-plane (which contains $\cal C$)  by {\it{horizontal plane}},
and any plane that contains the $z$-axis by {\it{vertical plane}}.\\

A {\it symmetry} $\psi$ of our problem is a transformation
$\psi :\R^3-{\cal C}\rightarrow\R^3-{\cal C}$
that preserve $V$, i.e $V$ is invariant by $\psi$, that is $V\psi =V$.
 The symmetries of $V$
are evident from the geometry of the problem: they are the rigid
transformations of $\R^{3}$ that preserve ${\cal C}$:

{\prop The potential $V$ is invariant by:\\
(i) Rotations about the z-axis,\\
(ii) Reflection with respect to the horizontal plane,\\
(iii) Reflections with respect to a vertical plane. \label{1.4.1}}
\vspace{0,2cm}

\noindent  {\bf  Proof.} It follows from (\ref{1.3}) and from the
expressions of $D$ and $d$ (see section 1.1). \CaixaPreta

\vspace{0,3cm}

{\cory The gradient field $\nabla V$ is invariant by:\\
(i) Rotations about the z-axis,\\
(ii) Reflection with respect to the horizontal plane,\\
(iii) Reflections with respect to a vertical plane.  \label{1.4.2}}

\vspace{0,2cm}

\noindent  {\bf  Proof.} Since all symmetries above are isometries
of $\R^{3}\setminus{\cal{C}}$ (with the canonical flat metric),
 the Corollary follows from the Proposition above and from a classical result in
the elementary theory of mechanical systems on Riemannian
manifolds.
 \CaixaPreta

\vspace{0,2cm}

{\cory If $\mbox\rv(t)$ is a solution of (\ref{00}), then
$\psi\,\mbox\rv (t)$ is also a solution of (\ref{00}), where
$\psi$ is one of the symmetries above. \label{1.4.3}\CaixaPreta}

\vspace{0,2cm}

{\cory  The z-axis, the horizontal plane and the vertical planes are
invariant subspaces, that is, a solution of (\ref{00})
that starts in one of these spaces, with velocity contained in it,
stays there for all the time in which the solution is defined.
\label{1.4.4}} \vspace{0,2cm}

 \noindent {\bf  Proof.}
 It follows from the following facts:

- the $z$-axis is the invariant space of all the symmetries in $(i)$
of Corollary \ref{1.4.2},

- the horizontal plane is the invariant space of the reflection in
$(ii)$ of Corollary \ref{1.4.2},

- a vertical plane is the invariant space of a refection in $(iii)$
of Corollary \ref{1.4.2}.
 \CaixaPreta\\

 Note that a {\it radial line}
(i.e. any line in the $xy$-plane that passes through the origin) is also an invariant subspace
 since it is the intersection of the horizontal plane with a vertical plane.

\vspace{0,5cm}

\subsection{Properties.}

We will need to consider circles with variable radius and mass.
Write  $V(\mbox{\rv},\rho,M)$ to denote
the potential induced by $\cal C$, contained in
the $xy$-plane and  centered at the origin, with radius $\rho$ and
mass $M$, and $\nabla
V(\mbox{\rv},\rho,M)$ to denote the gradient (with respect to
$\mbox{\rv}$) of $V(\mbox{\rv},\rho,M)$.

{\lem The   potential $V$ of the  fixed homogeneous
circle problem satisfies the following identities:

\noindent (i) $V(\mbox{\rv},\rho,c\,M)= c
\,V(\mbox{\rv},\rho,M),\,\,$ for  $c\in\R$,

\noindent (ii) $V(c\,\mbox{\rv},c\,\rho,M)= \frac{1}{c} \,
V(\mbox{\rv},\rho,M),\,\,$ for  $c>0$,

\noindent (iii) $\nabla V(c\,\mbox{\rv},c\,\rho,M)=
\frac{1}{c^{2}}\, \nabla V(\mbox{\rv},\rho,M),\,\,$ for
$c>0$,

\noindent (vi) $\nabla V(\mbox{\rv},\rho,c\,M)= c \,\nabla
V(\mbox{\rv},\rho,M),\,\,$ for  $c\in\R$. \label{1.3.1}}

\vspace{0,2cm}

\noindent  {\bf  Proof.} It follows directly from the definition
of $V(\mbox{\rv},\rho,M)=-\frac{M}{2\pi}
\int_{0}^{2\pi}\frac{d\theta}{|\!|\mbox{\rv} -\mbox{$\rho
e^{i\theta}$}|\!|}$, where $e^{i\theta}=(cos\theta,sin\theta,0)$.
\CaixaPreta

{\cory Let  $\rho,\zeta, M, N$ be positive numbers.
If  ${\mbox{\rv}}(t)$ is a solution of
$\,\stackrel{..}{\mbox{\rv}}(t)=-\nabla V({\mbox{\rv}},\rho,M)$
then  ${\mbox{\sv}}(t)=\frac{\zeta}{\rho}\mbox{\rv} \left(
\sqrt{\frac{N\rho^{3}}{M\zeta^{3}}}\,\,t \right)$ is a solution
of  $\,\stackrel{..}{\mbox{\sv}}(t)=-\nabla
V(\mbox{\sv},\zeta,N)$. \label{1.3.2}}

\vspace{0,2cm}

\noindent  {\bf  Proof.} It follows from Lemma \ref{1.3.1}, by a
direct calculation. \CaixaPreta

 {\obs {\rm   Note that if  ${\mbox{\rv}}(t)$ and
${\mbox{\sv}}(t)$ are as above, then they have the same qualitative
properties. For instance, if  ${\mbox{\rv}}(t)$ is periodic
 then ${\mbox{\sv}}(t)$ is also periodic (with
period $\frac{T\sqrt{M\zeta^{3}}}{\sqrt{N\rho^{3}}}$, where $T$ is the period of
${\mbox{\rv}}(t)$).}
\label{1.3.3}}\\

Note that Lemma \ref{1.3.1} and Corollary  \ref{1.3.2} imply
that in the study the fixed homogeneous circle problem we can assume the
mass and the radius to be equal to one.

\vspace{0,3cm}

\section{Singularities: Proof of Theorem A.}

To prove Theorem A we first show the following Proposition, which
is of general nature. In the next Proposition {\it dist} denotes
``Euclidean distance".

{\prop  Let $V:\Omega\rightarrow\R$, $\Omega\subseteq\R^{n}$ open.
Let ${\mbox{{\mbox{\rv}}}}(t)$, $t\in(a,b)$, be a solution of
$\stackrel{..}{\mbox{\rv}}\,\,=-\nabla V({\mbox{\rv}}) $ such that
there exist:

(1) $v_{0},v_{1}\in\R, \,\,v_{0}<v_{1}$, with
dist$\,(V^{-1}(v_{0}),V^{-1}(v_{1}))>0,$

(2) $ t_{1}<s_{1}<t_{2}<s_{2}\ldots,\,\,\,\,\,
t_{i},s_{i}\in(a,b),\,\,\, i\in\N$, satisfying the following two properties:
(a) $ V({\mbox{{\mbox{\rv}}}}(t_{i}))=v_{0}$,  and (b) $V({\mbox{{\mbox{\rv}}}}(s_{i}))=v_{1}.$
 Then $b=+\infty$. \label{2.0.8}} \vspace{0,3cm}

\noindent  {\bf  Proof.} Since $V$ is continuous, for every
interval $[t_{i},s_{i}]$ we can choose, for all  $i$,
 an interval
$[t_{i}^{\ast},s_{i}^{\ast}]\subseteq [t_{i},s_{i}]$ such  that
$v_{0}=V({\mbox{\rv}}(t_{i}^{\ast}))\leq V({\mbox{\rv}}(t))\leq
V({\mbox{\rv}}(s_{i}^{\ast}))=
v_{1},\,\,\mbox{for all}\,\, t\in[t_{i}^{\ast},s_{i}^{\ast}]$.

Since the  total energy $E = E({\mbox{{\mbox{\rv}}}}(t)) =
V({\mbox{{\mbox{\rv}}}}(t))+\frac{1}{2}|\!|\dot {\mbox{\rv}}
(t)|\!|^{2}$ is constant, we have that $|\!|\dot
{\mbox{\rv}}(t)|\!|= \sqrt{2(E-V)},$ with
$V=V({\mbox{{\mbox{\rv}}}}(t))$ and  $E$ constant. Then we have $|\!|\dot {\mbox{\rv}}(t)|\!|\leq
\sqrt{2(E-v_{0})},\,$ for $ t\in [t_{i}^{\ast},s_{i}^{\ast}]$.
Since the length of the curve ${\mbox{{\mbox{\rv}}}}(t)$ between
$[t_{i}^{\ast},s_{i}^{\ast}]$ is larger than the distance\\
 $d=
dist(V^{-1}(v_{0}),V^{-1}(v_{1}))>0$, we have
$$d\leq  \int_{t_{i}^{\ast}}^{s_{i}^{\ast}} |\!|\dot{{\mbox{\rv}}}(t)|\!|dt\,\leq
\int_{t_{i}^{\ast}}^{s_{i}^{\ast}}\sqrt{2(E-v_{0})}dt \,=
(s_{i}^{\ast}-t_{i}^{\ast})\sqrt{2(E-v_{0})}\leq
(s_{i}-t_{i})\sqrt{2(E-v_{0})}.$$

\noindent Hence, $\frac{d}{\sqrt{2(E-v_{0})}}\leq (s_{i}-t_{i}),\,\,\,
\mbox{for all}\,\,i.$
Since we have an  infinite number of disjoint intervals
$(t_{i},s_{i})$, with $0<\frac{d}{\sqrt{2(E-v_{0})}}
\leq(s_{i}-t_{i})$ for all $i$, we conclude that
$b=+\infty.$ \CaixaPreta \\

In what follows of this section
let $\cal C$ be the fixed homogeneous circle in $\R^{3}$
centered at the origin, contained in the horizontal plane and
with constant density $\lambda$ and radius 1. Also let $V$ be the potential
 induced by $\cal C$ and ${\mbox{\rv}}(t)$ be a solution of
$\stackrel{..}{\mbox{\rv}}\,=-\nabla V({\mbox{\rv}}),$ defined in the
maximal interval  $(a,b)$.

{\prop  If $\,\, b<+\infty$ then
$lim_{t\rightarrow b^{-}} dist({\mbox{\rv}}(t),{\cal{C}})=0$.
 Analogously, if $\,\, a>-\infty$ then
$lim_{t\rightarrow a^{+}} dist({\mbox{\rv}}(t),{\cal{C}})=0$.
\label{2.0.9}}\\

\noindent To  prove  this  Proposition we need the  following Lemmas.
Let $\{ {\mbox{\rv}}_{n}\}_{n\in\N} \subset \R^3-{\cal C}$.

{\lem  $lim_{n\rightarrow +\infty}V({\mbox{\rv}}_{n})= 0\,\,$ if
and only if $\,\,lim_{n\rightarrow
+\infty}|\!|{\mbox{\rv}}_{n}|\!|\,= +\infty$. \label{2.0.10}}
\vspace{.1in}

\noindent  {\bf  Proof.} Suppose first that
$|\!|{\mbox{\rv}}_{n}|\!|\rightarrow +\infty$. For all $u\in {\cal
{C}},$
 we have $|\!|u|\!|=\rho.$ Hence, for
$|\!|{\mbox{\rv}}|\!|> \rho$, and $u\in {\cal {C}},$ we have
$|\!|{\mbox{\rv}}-u|\!|\,\geq | |\!|{\mbox{\rv}}
|\!|-|\!|u|\!||\,= ||\!|{\mbox{\rv}}|\!|-\rho| \,\geq
|\!|{\mbox{\rv}}|\!|-\rho$.  Hence
$0\leq -V({\mbox{\rv}})=\int_{{\cal
{C}}}\frac{\lambda}{|\!|{\mbox{\rv}}-u|\!|}du\, \leq \int_{{\cal
{C}}}\frac{\lambda}{|\!|{\mbox{\rv}}|\!|-\rho}du\,=
\frac{M}{|\!|{\mbox{\rv}}|\!|-\rho},$
where $M=\int_{{\cal {C}}}\lambda \,du\,$ is the mass of ${\cal {C}}$.
It follows that $lim_{n\rightarrow +\infty}V({\mbox{\rv}}_{n}) =0.$

Suppose now that $V({\mbox{\rv}}_{n})\rightarrow 0$. Since
$|\!|u|\!|=\rho$, for  $u\in {\cal {C}}$, we have
$\,\frac{1}{|\!|{\mbox{\rv}}-{u}|\!|}\geq
\frac{1}{|\!|{\mbox{\rv}}|\!|+\rho }\,\,$. Consequently,
$0\leq \frac{M}{|\!|{\mbox{\rv}}|\!|+ \rho}=\int_{{\cal {C}}}\frac{\lambda }
{|\!|{\mbox{\rv}}|\!|+\rho}du\,\leq \int_{{\cal
{C}}}\frac{\lambda}{|\!|{\mbox{\rv}}-u|\!|}du\,=
-V({\mbox{\rv}}).$
Therefore, $lim_{n\rightarrow
+\infty}|\!|{\mbox{\rv}}_{n}|\!|=+\infty$. This proves the Lemma. \CaixaPreta

\vspace{0,3cm}

{\lem   $lim_{n\rightarrow +\infty} V({\mbox{{\mbox{\rv}}}}_n)=
-\infty\,\,$ if and only if  $\,\,lim_{n\rightarrow +\infty}
dist({\mbox{{\mbox{\rv}}}}_n,{\cal {C}})=0$. \label{2.0.11}}
\vspace{0,2cm}

\noindent  {\bf  Proof.} Suppose that $lim_{n\rightarrow +\infty}
V({\mbox{{\mbox{\rv}}}}_n)= -\infty$.  Since $0<
dist({\mbox{{\mbox{\rv}}}}_n,{\cal {C}})\leq |\!|{\mbox{\rv}}_{n}
-u|\!|,\,\,\,\mbox{for all}\,\,\, u\in {\cal {C}},$ we have $0<
-V({\mbox{{\mbox{\rv}}}}_{n})=\int_{{\cal {C}}}\frac{\lambda}
{|\!|{\mbox{\rv}}_{n} -u|\!|}du\,\leq \int_{{\cal
{C}}}\frac{\lambda}{dist({\mbox{\rv}}_{n},{\cal {C}})}du\,=
M\frac{1}{dist({\mbox{\rv}}_{n},{\cal {C}})}.$ It follows that
$dist({\mbox{{\mbox{\rv}}}}_n,{\cal {C}})\rightarrow 0$.

Conversely,  suppose  that    $dist({\mbox{{\mbox{\rv}}}}_n,{\cal
{C}})\rightarrow 0$. Let $d, D$ be as in section 1.1. We have
 $d({\mbox{{\mbox{\rv}}}}_n)=
dist({\mbox{{\mbox{\rv}}}}_n,{\cal {C}})\rightarrow 0$ and
$D({\mbox{\rv}}_n)\rightarrow 2\rho$. It follows from   ($vi$) of
Lemma \ref{1.6.2} that
$V({\mbox{\rv}}_n)=\frac{-4\lambda}{D({\mbox{\rv}}_n)}\,f\bigl(
\frac{d^2 ({\mbox{\rv}}_n)}{D^2 ({\mbox{\rv}}_n)}\bigr)$
$\rightarrow -\infty.$ \CaixaPreta

\vspace{0,6cm}

\noindent {\bf Proof of  Proposition \ref{2.0.9}.} First note
that, by the two Lemmas above, $V^{-1}(c)$ is compact, for all
$c<0.$ Now, let   ${\mbox{\rv}}(t),\,\, t\in (a,b),\,\,
b<+\infty$, be a maximal solution of
$\stackrel{..}{\mbox{\rv}}\,=-\nabla V({\mbox{\rv}})$. Then
$lim_{t\rightarrow b^{-}}V({\mbox{{\mbox{\rv}}}}(t))$ exists (it
could be finite or infinite), otherwise, by the continuity of  $V$
we could choose two sequences $(t_{n})$, $(s_{n})$ such that
${t}_{1}<{s}_{1}<{t}_{2}<{s}_{2}<\ldots,\,\,{t}_{i},{s}_{i}\in(a,b),$
with $V({\mbox{{\mbox{\rv}}}}(t_{i}))=v_{0},
V({\mbox{{\mbox{\rv}}}}(s_{i}))=v_{1},\,v_{0}<v_{1}$. Since
$V^{-1}(v_{0})$ and $V^{-1}(v_{1})$ are   disjoint compact and
non-empty, we have that $dist(V^{-1}(v_{0}),V^{-1}(v_{1}))>0.$
Hence, by Proposition  \ref{2.0.8} we would have that
$b=+\infty$, a contradiction.

Since $image\, V\subset (-\infty ,0)$, it follows that
$lim_{t\rightarrow b^{-}}V({\mbox{{\mbox{\rv}}}}(t))$  is either
$0$, or  it is a number $v^{\ast}\neq 0$, or it is equal to
$-\infty$. We show  that the two first possibilities do not
happen. It will then follow that $lim_{t\rightarrow
b^{-}}V({\mbox{{\mbox{\rv}}}}(t))=-\infty$ and therefore, by Lemma
\ref{2.0.11}, we would have that $lim_{t\rightarrow
b^{-}}dist(\mbox{\rv}(t),{\cal C})=0,$ which proves the
Proposition.

If $lim_{t\rightarrow b^{-}}V({\mbox{{\mbox{\rv}}}}(t))=0,$
then there exists  $t_{0}$ such  that for all
$t>t_{0},\,\,v_{0}<V({\mbox{{\mbox{\rv}}}}(t))<0,$ with $v_0
=V({\mbox{{\mbox{\rv}}}}(t_{0})).$ On the other hand by Lemma \ref{2.0.10}, we have
that $|\!|{\mbox{\rv}}(t)|\!|\rightarrow +\infty$ when
$t\rightarrow b^{-}.$  Moreover, since
$E=V({\mbox{{\mbox{\rv}}}}(t))+\frac{1}{2}|\!|\dot{\mbox{\rv}}(t)|\!|^{2},\,\,E $
is constant
then $|\!|\dot{\mbox{\rv}}(t)|\!|<
\sqrt{2(E-v_{0})},\,\,\mbox{for all}\,\, t\in (t_{0},b),$ that is,
the  velocity is bounded in this interval.
We also have that
$\int_{t_{0}}^{t}|\!|\dot{{\mbox{{\mbox{\rv}}}}}(t)|\!|dt\,\geq
|\!|{\mbox{{\mbox{\rv}}}}(t)-{\mbox{{\mbox{\rv}}}}(t_{0})|\!|\,\geq
dist({\mbox{{\mbox{\rv}}}}(t), V^{-1}(v_{0}))=: d_{t}.$
Hence, $d_{t}\leq (t-t_{0})\sqrt{2(E-v_{0})},$ and  then,
$ t_{0}+\frac{d_{t}}{\sqrt{2(E-v_{0})}}\leq t,
\,\,\,\mbox{for all}\,\,\,t\in(t_{0},b).$

Since $V^{-1}(v_{0})$ is compact, there exists  $r>0$ such  that
$V^{-1}(v_{0})\subset B(0,r)$ and  since
$|\!|{\mbox{\rv}}(t)|\!|\rightarrow +\infty,$ given  $n>0,$ there exists
$t_{n}$ such that for all $t\in[t_n,b)$,
${\mbox{{\mbox{\rv}}}}(t)\notin B(0,n).$ This implies
that, $d_{t}=dist({\mbox{{\mbox{\rv}}}}(t),V^{-1}(v_{0}))>n-r,$
 for all
$b>t\geq t_{n}$. In particular, $d_{t_{n}}\geq n-r$, and for
$b>t_{n}>t_{0}$, we have $b\geq t_{n}\geq
t_{0}+\frac{n-r}{\sqrt{2(E-v_{0})}}.$
Hence, $lim_{n\rightarrow \infty}\left(
t_{0}+\frac{n-r}{\sqrt{2(E-v_{0})}}\right)=+\infty,$
 and we conclude that $b=+\infty$, a contradiction. \vspace{0,2cm}

Finally suppose that $lim_{t\rightarrow
b^{-}}V({\mbox{{\mbox{\rv}}}}(t))=v^{\ast},$ with
$v^{\ast}\in(-\infty,0)$. Hence, given  $\varepsilon>0,$ there exists
$t_{\varepsilon}$ such  that for all
$t>t_{\varepsilon},\,\,V({\mbox{{\mbox{\rv}}}}(t))\in(v^{\ast}-\varepsilon,
v^{\ast}+\varepsilon).$ Let  $v_{0}=v^{\ast}-\varepsilon
,\,\,v_{1}=v^{\ast}+\varepsilon;$ we can suppose that $0<\varepsilon
<|v^{\ast}|,$ hence $v_{1}< 0.$ Clearly the set
$V^{-1}([v_{0},v_{1}])$ is not empty and since $V$ is continuous,
 this set  is closed in $\R^{3}\setminus \cal C$. Moreover,
by Lemma \ref{2.0.11}, this set is closed in $\R^{3}$, and by
Lemma \ref{2.0.10} it is bounded (because $v_{1}<0$). Hence,
$V^{-1}([v_{0},v_{1}])$ is compact and
${\mbox{{\mbox{\rv}}}}(t)\in V^{-1}([v_{0},v_{1}]),\,\,\mbox{ for
all}\,\,\,t>t_{\varepsilon}.$

On the other hand, from  the energy equation, we have
$|\!|\dot{{\mbox{{\mbox{\rv}}}}}(t)|\!|^{2}=
2(E-V({\mbox{{\mbox{\rv}}}}(t)))$, therefore,   for all
$t>t_{\epsilon}$, we have $\,|\!|\dot{{\mbox{\rv}}}(t)|\!|\leq
\sqrt{2(E-v_{0})}= c_{2}.$
Hence the  maximal solution
$({\mbox{{\mbox{\rv}}}}(t),\dot{{\mbox{{\mbox{\rv}}}}}(t)),$
$\,t\in(a,b),$ of the  system of first order differential
equations \vspace{0,15cm}
$$ \left\{\begin{array}{l}
\dot{\mbox{\rv}}\,\,=\,\,{\mbox{\vv}}\\
\dot {\mbox{\vv}}\,\, =\,\, -\nabla V({\mbox{\rv}})
\end{array}\right.$$

\noindent is contained in  the compact
$V^{-1}([v_{0},v_{1}])\times \overline {B(0,c_{2})},$ for $t\in(t_{\epsilon}, b).$
It follows from a classical result in the elementary
theory of differential equations (see \cite{S}) that $b=+\infty,$
a contradiction. This proves the Proposition. \CaixaPreta \\

\noindent  {\bf  Proof of Theorem A.} Let ${\mbox{\rv}}(t)$,
$t\in(a,b)$, $b<+\infty$, be  a maximal solution of
 $\stackrel{..}{\mbox{\rv}}\,\,=\,-\nabla
V({\mbox{\rv}})$. By Proposition \ref{2.0.9}, we have that
$lim_{t\rightarrow b^{-}} dist({\mbox{{\mbox{\rv}}}}(t),{\cal{C}}
)=0$.
 To  prove  that $lim_{t\rightarrow b^{-}} {\mbox{\rv}}(t)={\mbox{\rv}}^{\ast}\in
\cal C$, we write this system in cylindrical coordinates
$(r,\theta, z)$ in ($\R^{3}\setminus\{z-\mbox{axis}\}$). We have
(see equations (\ref{R4}) of section 5):
\vspace{0,1cm}
$$\stackrel{..}r\,\,=\,\, \frac{K^{2}}{r^{3}}-\frac{\partial V}{\partial
r},\,\,\,\,\,\,\,\,\, \stackrel{..}z\,\, =\,\, -\frac{\partial
V}{\partial z},$$

\noindent with $\dot\theta =\displaystyle\frac{K}{r(t)^{2}}$, where
$K$ is constant and $V(r,z)=V(r\,cos\theta,\,r\,sin\theta,z)$.

Since the circle in cylindrical coordinates is given  by ${\cal{
C}} =\{ (\rho,\varphi,0), \varphi\in\R\},$
 we have that showing  $lim_{t\rightarrow b^{-}} {\mbox{{\mbox{\rv}}}}(t)=
 {\mbox{{\mbox{\rv}}}}^{\ast}\in \cal C$ is
equivalent to showing $lim_{t\rightarrow b^{-}} r(t)=\rho$,
$lim_{t\rightarrow b^{-}} z(t)=0$ and $lim_{t\rightarrow b^{-}}
\theta(t)=\theta_{0}$, for some $\theta_{0}$. The two first limits
follow from the fact that $lim_{t\rightarrow b^{-}}
dist({\mbox{\rv}}(t),{\cal C})=0$. We will now  prove  that
$lim_{t\rightarrow b^{-}} \theta(t)=\theta_{0} $.
If $K=0$, $\theta(t)$ is constant, and we have nothing to  prove.
Suppose then that $K>0$. Hence $\theta(t)$ is an increasing function.
 Thus, to prove that the limit
of $\theta(t)$  exists, it is enough to prove  that $\theta(t)$ is
bounded above, for $t$ in a neighborhood of
 $b$.

 Since $r(t)\rightarrow \rho,$ there exists  $t_{0}$ such  that, for
$t>t_{0},\,\,r(t)>\frac{\rho}{2}.$ Hence
$ \theta(t)
=  \displaystyle\int_{t_{0}}^{t}\frac{K}{r(s)^{2}}ds +
\theta(t_{0})
\leq   \displaystyle\int_{t_{0}}^{t}\frac{4 K}{\rho^{2}}ds +
\theta(t_{0})
\leq   \displaystyle\frac{4 K (b-t_{0})}{\rho^{2}}+
\theta(t_{0}) <+\infty $
for all $t\in(t_{0},b)$.
 Therefore, the limit of $ \theta(t)$  when   $t\rightarrow
b^{-}$  exists. This proves the Theorem.
\CaixaPreta
\vspace{.4in}

\section {The Dynamics Inside the Circle: Proof of Theorem B.}

In this section we consider again the fixed homogeneous circle
$\cal C$  contained in the $xy$-plane, centered at the origin and
with radius 1 and constant density $\lambda$. Let
${\mbox{\rv}}=(x,y)$ be the position of the particle  in this
plane, under the influence of the gravitational attraction induced
by $\cal C$. Let  $r\,=|\!|\mbox{\rv} |\!|$ and $\theta$ be the
polar coordinates of ${\mbox{\rv}}=(x,y)$. Also, let $D$ and $d$
be as in section 1.1. In the horizontal plane we have
$D^{2}=(r+1)^{2}$ and $d^{2}=(r-1)^{2}=(1-r)^2$. Hence the expression
for the potential (\ref{1.3}) becomes
$$V({\mbox{\rv}})=V(r)=-4\lambda\displaystyle{\int_{0}^{\frac{\pi}{2}}} \frac{d\theta}
{\sqrt{({r} +1)^{2}cos^{2}\theta + ({r} -1)^{2}sin^{2}\theta}}$$

Since $V$ depends only on $r$, we have a central force problem and (see \cite{G}, ch.3)
the system of equations
$\stackrel{..}{\mbox{\rv}}\, =-\nabla
V({\mbox{\rv}})$ (restricted to the horizontal plane) is
equivalent to the system:

\begin{equation} \left\{ \begin{array}{l}
\stackrel{..}{{r}}\,\,=\,\, \displaystyle\frac{K^{2}}{{r}^{3}}-\frac{d}{dr}V(r)\\
\\
\stackrel{.}\theta\,\, =\,\, \displaystyle\frac{K}{{r}^2}\end{array}
\right.
\label{6}
\end{equation}

\noindent where $K=\langle\,{\mbox{\rv}} \times
\dot{\mbox{\rv}},\,e_{3}\,\rangle =\,x{\dot y} -{\dot x} y$
is the angular momentum. Note that  $\stackrel{..}{{r}} \,=\frac{K^{2}}{{r}^{3}}-\frac{d}{dr}V(r)=
\frac{d}{d{{r}}}\left( \frac{-K^{2}}{2{r}^{2}}-V({r})\right)$ which is equivalent to
\begin{equation}
\stackrel{..}{{r}} \,=-\frac{d}{dr}U({r}),
\label{4.0}
\end{equation}
 where
$U({r})=\displaystyle\frac{K^{2}}{2{r}^{2}}+V({r})$. The function $U(r)$ is called
the {\it effective potential }.\\

{\prop    For $0<r<1,$
$\frac{d}{dr}V({r})<0$. \label{4.2.2}}\\

\noindent  {\bf  Proof.} In the horizontal plane and inside the
circle we have $D=r+1$ and $d=1-r$. Hence $D_1=1$ and
$d_1=\sqrt{1-r^2}$, where $D_1$ and $d_1$ are the arithmetic and
geometric means of $D$ and $d$, respectively. Therefore, by Remark
\ref{1.1} of section 1.1, we have
\begin{equation}
V({\mbox{\rv}})=V(r)=-4\lambda T(1+r,1-r)= -4\lambda
T(1,\sqrt{1-r^2})=-4\lambda\displaystyle{\int_{0}^{\frac{\pi}{2}}}
\frac{d\theta} {\sqrt{1-r^{2}sin^{2}\theta}}.
\label{3.6}\end{equation}

\noindent Differentiating we have
$\frac{d}{dr}V(r)=-4\lambda\displaystyle{\int_{0}^{\frac{\pi}{2}}}
\frac{\, r\, sin^2 \theta \,\,\,\,\,}{\Bigl( 1-r^{2}sin^{2}\theta\Bigr)^{3/2}  }\,
 d\theta\,\, <\,\, 0.$
 \CaixaPreta
\vspace{0,4cm}

\noindent   {\bf Proof of Theorem B.}
Differentiating we have $\frac{d}{dr}U(r)=-\frac{K^2}{r^3}+\frac{d}{dr}V(r)$.
Hence, by the Proposition above we have  $\frac{d}{dr}U(r)<0$ for $0<r<1$.
This proves ($i$) of Theorem B.

When $r\rightarrow 1^-$ we have that the particle tends to the circle. Hence,
by Lemma \ref{2.0.11}, $V(r)$ tends to $-\infty$. Consequently
{\footnotesize $U({r})=\displaystyle\frac{K^{2}}{2{r}^{2}}+V({r})$} tends also
to $-\infty$.
This proves ($ii$). Note that ($iii$) and ($iv$) follow from
 the definition of $U(r)$ and (\ref{3.6}). \CaixaPreta
\vspace{.4in}

\section{Proofs of Propositions 1 and 2.}

First we prove Proposition 2, but before  we give some remarks.
Let ${\mbox{\rv}}(t)$ be as  in the statement of Proposition 2.
Without loss of generality we can assume that $t_0=0$
(see definition of a normalized solution).\\

\noindent (1) Recall that the curvature of the curve ${{\mbox{\rv}}}(t)$ is
by definition
$$k(t)=\frac{1}{|\!|\dot{{\mbox{\rv}}}(t)|\!|^{2}}
\langle{\stackrel{..}{\mbox{\rv}}}(t),{\mbox\nn} (t)\,\rangle\,=\,
\frac{1}{|\!|\dot{{\mbox{\rv}}}(t)|\!|^{3}}
\,\langle\,{\stackrel{..}{\mbox{\rv}}}(t),
\,|\!|\dot{{\mbox{\rv}}}(t)|\!|\,{\mbox\nn} (t)\rangle\,=\,
\frac{1}{|\!|\dot{{\mbox{\rv}}}(t)|\!|^{3}}
\,\langle\,{\stackrel{..}{\mbox{\rv}}}(t),{\dot{\mbox{\rv}}}^{\perp}(t)\,\rangle,$$

\noindent
where ${\mbox\nn} (t)=\frac{(-\dot y (t),\dot x (t))}{|\!|(-\dot y (t),\dot x (t))
|\!|}$
and
$\,\,{\dot{\mbox{\rv}}}^{\perp}(t)=\|{\dot{\mbox{\rv}}}(t)\|\,{\mbox\nn} (t)=\,(-\dot
y (t),\dot x(t)).$\\

\noindent (2)  Note that $k(t)$ is continuous and is never zero
because if $k(t)=0$ for some $t$, then
${\stackrel{..}{\mbox{\rv}}}(t) \bot\, {\mbox\nn} (t)$ and follows
that $\dot{\mbox{\rv}}(t)$ is radial (because
${\stackrel{..}{\mbox{\rv}}}(t)$ is always radial), and then the
solution  ${\mbox{\rv}}(t)$
would be radial, that is  $K=0$, a contradiction.\\

\noindent (3) $\dot{{\mbox{\rv}}}(0)\,\bot
{\stackrel{..}{\mbox{\rv}}} (0)$ and
${\stackrel{..}{\mbox{\rv}}}(0)$ points upward. It follows that
${\mbox\nn} (0)=c{\stackrel{..}{\mbox{\rv}}} (0),\,$ for some $c>0$.
Then $k(0)=\bar{c}\,
 |\!|{\stackrel{..}{\mbox{\rv}}}(0)|\!|^{2}>0$, with $\bar{c}>0$. In this
 case, since  $k$ is continuous and is always non-zero, we have that
 $k>0$, for all $t$.\\

\noindent {\bf Proof of  Proposition 2. } Let ${\mbox{\rv}} (t)=
(x(t), y(t))$ be as in the statement  of Proposition 2 and recall
that we are assuming $t_0=0$. Let $(a,b)$ be the maximal interval
on which
$\mbox{\rv}(t)$ is defined. Then $a<0<b$.\\

\noindent {\bf Claim.} $\dot{x} (t)\neq 0,$ for all $t\in (a,b)$.

We prove that $\dot{x} (t)\neq 0,$ for all $t\in [0,b)$, the proof
for $t\in (a,0]$ is similar. Suppose that there exists  $t_0 \in
[0,b)$ such  that $\dot{x}(t_{0})=0.$ Let $t_{0}=min\{t\geq 0\,
;\,\dot{x}(t)=0\}.$ Since $\dot{x}$ is continuous and
$\dot{x}(0)>0$, it follows that $t_{0}>0.$ Then $\dot{x}(t)>0,$ for
$t\in[0,t_{0}).$ We have
$\dot{{\mbox{\rv}}}(t_{0})=(\dot{x}(t_{0}),\dot{y}(t_{0}))=(0,\dot{y}(t_{0}))$,
with $\dot{y}(t_{0})\neq 0,$ because $\dot{{\mbox{\rv}}}(t)\neq 0$,
for all $t$. We can write $\dot{{\mbox{\rv}}}(t)=a(t)(cos\varphi(t),
sin\varphi(t)),$ with $a(t)= \|\dot{{\mbox{\rv}}}(t)\|>0$ and
$\varphi(0)=0$, $\varphi$ continuous. Differentiating:
${\stackrel{..}{\mbox{\rv}}}
(t)=\dot{a}(t)(cos\varphi(t),sin\varphi(t))+
a(t)\dot{\varphi}(t)(-sin\varphi(t),cos\varphi(t))$. Hence
$k(t)=\,\frac{1}{|\!|\dot{{\mbox{\rv}}}(t)|\!|^{3}}
\langle\,{\stackrel{..}{\mbox{\rv}}}(t),
\dot{{\mbox{\rv}}}^{\bot}(t)\,\rangle \,=\frac{a^{2}(t)\dot\varphi
(t)}{a^{3}(t)}\,=\,\frac{\dot\varphi (t)}{a(t)}.$ Because the
curvature is positive (see Remark (3) above) this shows that
$\dot{\varphi}(t)>0$ and $\varphi(t)$ is an increasing function.
Since $\dot{x}(t_{0})=0$, $\,\varphi(t_{0})=\frac{\pi}{2}$ (hence
$\dot{y}(t_{0})>0$), or $\,\varphi(t_{0})=\frac{3\pi}{2}$ (hence
$\dot{y}(t_{0})<0$). Since $\varphi$ is an increasing function we
have that $\varphi(t_{0})=\frac{\pi}{2}$, which implies that
$\dot{y}(t_{0})>0$. Since  $x$ is increasing on  $(0,t_{0})$, we
have that
 $x(t_{0})>0.$
Recall that   ${\stackrel{..}{\mbox{\rv}}}(t)$ is radial and
expansive, that is,
$\,{\stackrel{..}{\mbox{\rv}}}(t)=b(t){\mbox\rv}(t)
=b(t)(x(t),y(t)),$ with $b(t)>0$. Hence $k(t_{0})=b(t_{0})
\,\langle\,(x(t_{0}),y(t_{0})),
(-\dot{y}(t_{0}),0)\,\rangle\,=\,b(t_{0})(-x(t_{0})\dot{y}(t_{0}))<0,$
a  contradiction. Therefore   there is no $t_{0}$ such  that
$\dot{x}(t_{0})=0$. This proves the claim.  \CaixaPreta
\vspace{0,3cm}

It follows from the claim that $\dot{x}(t)>0$, $\,a<t<b.$
Hence $x(t)$  is a increasing function.
In this way the function $t\rightarrow x(t)$ is one-to-one and it
follows that $x(t)$ possesses an inverse $t=t(x)$.
Define  $f(x)=y(t(x))$. Note that the graph  of $f$ is equal to the
trace of  ${\mbox{\rv}}$.
Differentiating $f$ with respect to $x$, we have
$$\frac{d}{dx}f(x)=\frac{d}{dt}y(t(x))\frac{d}{dx}t(x)=\frac{\frac{d}{dt}y(t(x))}
{\frac{d}{dt}x(t(x))}.$$

Differentiating again we have
$$\frac{d^{2}}{dx^{2}}f(x)=\frac{1}
{{\dot x}^{3}}\,\langle\,\stackrel{..}{\mbox{\rv}}
,|\!|\dot{\mbox{\rv}}|\!|\,{\mbox\nn} \,\rangle \,=
\frac{k|\!|\dot{\mbox{\rv}}|\!|^{3}}{{\dot x}^{3}}>0$$

\noindent because $\dot{x}>0$ and  $k>0$.
Hence the trace of ${\mbox{\rv}}(t)$ is given by the graph  of a
convex function $f$.
By the symmetry of the problem, the solution  ${\mbox{\rv}} (t)$ is
symmetric  with respect to the  $y$-axis.
Hence the function $f$ is even. \CaixaPreta

 \vspace{.3in}

\noindent {\bf Proof of Proposition 1.} Let ${\mbox{\rv}}(t)$ and
$(a,b)$ be as  in the statement of Proposition 1. Let $r(t)=||
{\mbox{\rv}}(t)||$. Then $r(t)$ satisfies
$\stackrel{..}{{r}} \,=-\frac{d}{dr}U({r})$. Suppose that
${\mbox{\rv}}(t)\rightarrow 1^{-}$, when   $t\rightarrow b^{-}$.
By Theorem B we have $lim_{t\rightarrow b^{-}}U(r(t))=lim_{r\rightarrow
1^{-}}U(r)=-\infty$. Thus there exists $t_{0}$ such  that
$U(r(t))<E-1,$ for $t\in[\,t_{0}, b)$, where
$E=\frac{1}{2}\dot r ^{2}+U(r)$  is the energy  of $r(t)$. Hence,
for $t\geq t_{0},$ $\,\dot r (t)\neq 0$. Moreover, $\dot r(t)>0$
(because $r(t)\rightarrow 1^{-}$). Consequently, $r(t)$ is
one-to-one on $[\,t_{0}, b)$, and has an inverse $t=t(r)$, $r_{0}
\leq r<1,$ $\,r_{0}= r(t_{0}).$ From $E=\frac{1}{2}\dot
r^{2}+U(r)$, we have that $\dot r =\sqrt{2(E-U(r))}$. Then
{
$\,\int_{t_{0}}^{t(r)}dt\,=\,\int_{r_{0}}^{r}\frac{dr}{\sqrt{2(E-U(r))}}$}
and it follows that {
$t(r)\,=\,\int_{r_{0}}^{r}\frac{dr}{\sqrt{2(E-U(r))}}+t_{0}$},
$r_{0}\leq r<1.$  In this way $b= lim_{r\rightarrow 1^{-}} t(r)=$
{ $ \int_{r_{0}}^{1}
 \frac{dr}{\sqrt{2(E-U(r))}}+t_{0}
 < \int_{r_{0}}^{1}\frac{dr}{\sqrt{2}}+t_{0}
 \,\leq \,\frac{1}{\sqrt{2}}(1-r_{0}) \,+\,t_{0}\,<\,+\infty.$}
This proves part (1) of Proposition 1.\\

We now prove part (2) of Proposition 1.
Again, let ${\mbox{\rv}}(t)=(x(t),y(t))$ and $(a,b)$ be as  in the statement of Proposition 1.
Assume  $\|\mbox{\rv} (t) \|\rightarrow 1,$
when   $t\rightarrow b^{-}$.
We have
$E=\frac{1}{2}|\!|\dot{\mbox{\rv}}|\!|^{2}+V(\mbox{\rv})$. Then
$|\!|\dot{\mbox{\rv}}|\!|^{2}=2(E-V(\mbox{\rv})).$
When  $t\rightarrow b^{-},$ $|\!|{\mbox{\rv}}|\!|\rightarrow 1$ hence
$V(\mbox{\rv})\rightarrow -\infty$.  Therefore
$|\!|\dot{\mbox{\rv}}|\!|\rightarrow +\infty$.
Without loss of generality suppose  that  ${\mbox{\rv}}(t)$ converges
to the point  $(1,0)$ of the circle, that is
$x(t)\rightarrow 1$ and $y(t)\rightarrow 0$ when   $t\rightarrow b^{-}.$

We shall show that $\dot{\mbox{\rv}}(t)$
becomes  horizontal as $t\rightarrow
b^{-},$ that is  $\frac{\dot y (t)}{\dot x (t)}\rightarrow 0$
when   $t\rightarrow b^{-}.$ First, note that $\dot x (t)\neq 0$
for $t$ close to $b$.
Moreover, $\dot x (t)\rightarrow +\infty$ when   $t\rightarrow
b^{-}.$
To see this suppose that there exists a
sequence $t_{n}\rightarrow b^{-}$ with $|\dot x (t_{n})|<M$
for some $M$.
Then $|\dot y (t_{n}) |\rightarrow +\infty$ (because $|\!| \dot
{\mbox{\rv}} (t_{n}) |\!|\rightarrow +\infty$).
Hence $|K|=\left|\dot x (t_{n}) y(t_{n}) -\dot y (t_{n})
x(t_{n})\right| \rightarrow +\infty,$ a contradiction because
$K$ is constant.

Now, suppose that  $lim_{t\rightarrow b^{-}}\frac{|\dot
y (t)|}{|\dot x (t)|}\neq 0$. This implies  that there exists a
sequence $\{ t_{n}\}$ with $t_{n}\rightarrow b^{-}$
such  that $\frac{|\dot y (t_{n})|}{|\dot x (t_{n})|}\geq \delta,$
for some $\delta>0.$
Hence, we have $|K| = $ {\footnotesize $ | x (t_{n})\dot y (t_{n}) -y (t_{n})\dot x (t_{n})|\,
\geq |x (t_{n})| |\dot y (t_{n})| -|y (t_{n})||\dot x (t_{n})|
\geq |x (t_{n})|(\delta|\dot x (t_{n})|)-|y (t_{n})||\dot x (t_{n})|
= |\dot x (t_{n})|\Bigl( \delta|x (t_{n})|-|y (t_{n})|\Bigr).$}
Taking the limit when   $t_{n}\rightarrow b^{-}$, we have that  $|K|\rightarrow
+\infty,$ a contradiction.  Therefore, $lim_{t\rightarrow b}\frac{|\dot
y (t)|}{|\dot x (t)|}=0$. \CaixaPreta
\vspace{0,4cm}

\section{Dynamics Outside the Circle: Proof of Theorem C.}

As in section  3 we consider the fixed homogeneous circle $\cal C$  contained
in the $xy$-plane, centered at the origin and with radius 1 and constant density $\lambda$.
Also ${\mbox{\rv}}$ will denote the position of a particle  in this
plane, under the influence of the gravitational attraction induced by $\cal C$.
Let  $r\,=|\!|\mbox{\rv} |\!|$. \\

\noindent {\bf Proof of Theorem C.} Let $D=D(\mbox \rv)$ and
$d=d(\mbox \rv)$ be as in section 1.1, that is they are the maximum
and minimum distances from $ \mbox \rv$ to the circle. In the
horizontal plane and outside the circle we have $D=r+1$ and $d=r-1$.
Hence $D_1=r$ and $d_1=\sqrt{r^2-1}$, where $D_1$ and $d_1$ are the
arithmetic and geometric means of $D$ and $d$, respectively.
Therefore, by Remark \ref{1.1} of section 1.1, we have

{\footnotesize $$V({\mbox{\rv}})=V(r)=-4\lambda \,T(r+1,r-1)=-4\lambda \,T(r,\sqrt{r^2-1})=
-4\lambda\displaystyle{\int_{0}^{\frac{\pi}{2}}} \frac{d\theta}
{\sqrt{r^{2}cos^{2}\theta + (r^2 -1)sin^{2}\theta}}=
-4\lambda\displaystyle{\int_{0}^{\frac{\pi}{2}}}
\frac{d\theta} {\sqrt{r^{2}-sin^{2}\theta}}$$}

\noindent Differentiating we have

{\footnotesize
$$\frac{d}{dr}V(r)=4\lambda\displaystyle{\int_{0}^{\frac{\pi}{2}}}
\frac{ r \,\,\,\,\,}{\Bigl( r^{2}-sin^{2}\theta\Bigr)^{3/2}  }\,
d\theta\,\, >\,\, 0.$$}

\noindent This proves ($i$) of Theorem C.
When $r\rightarrow 1^+$ we have that the particle tends to the circle. Hence,
by Lemma \ref{2.0.11}, $V(r)$ tends to $-\infty$. Consequently
{\footnotesize $U({r})=\displaystyle\frac{K^{2}}{2{r}^{2}}+V({r})$} tends also to $-\infty$.
This proves ($ii$). Note that ($iii$) follows from the definition of $U(r)$ and Lemma \ref{2.0.10}.\\



\noindent We now prove ($iv$)-($ix$) of the statement of Theorem C.
Since by definition $U(r)=\frac{K^2}{2r^2}+V(r)$ we have that
$\frac{d}{dr}U(r)=0$ if and only if
$\frac{K^2}{r^3}=\frac{d}{dr}V(r)$, or equivalently
$K^2=r^3\frac{d}{dr}V(r)$. Define $g(r)=r^3\frac{d}{dr}V(r)$, $\,
r>1$. By the formula above we have {\footnotesize $g(r)\, =\,
4\lambda\displaystyle{\int_{0}^{\frac{\pi}{2}}} \frac{ r^4
\,\,\,\,\,}{\Bigl( r^{2}-sin^{2}\theta\Bigr)^{3/2}  }\, d\theta.$}

{\lem   The function $g(r)$ has the following properties
\begin{enumerate}
\item[{a.}]   $lim_{r\rightarrow 1^{-}}\, g(r)=\infty$,

\item[{b.}]  $\lim_{r\rightarrow \infty} \frac{g(r)}{r}=2\pi
\lambda$,

\item[{c.}]  $g(r)$ has exactly one critical point $r_0$, $1<r_0<2$.
\end{enumerate} \label{5a}}

\vspace{.2in} \noindent {\bf Proof.} If $r\rightarrow 1^-$ the
particle approaches the circle therefore, by Lemma \ref{2.0.11}
$V(r)\rightarrow -\infty$. Hence $\frac{d}{dr}V(r)$ is not bounded
when $r\rightarrow 1^-$. But
$\frac{d^2}{dr^2}V(r)-4\lambda{\int_{0}^{\frac{\pi}{2}}} \frac{ 2r^2
+ sin^2\theta }{( r^{2}-sin^{2}\theta)^{5/2} }\,
 d\theta < 0$ therefore
$\frac{d}{dr}V(r)\rightarrow \infty$, as $r\rightarrow 1^-$. This
proves ($a$.). For ($b$.) we have  $\frac{g(r)}{r}=$ {\footnotesize
$4\lambda{\int_{0}^{\frac{\pi}{2}}} \frac{ r^3 \,\,\,\,\,}{(
r^{2}-sin^{2}\theta)^{3/2}  }\, d\theta=
4\lambda{\int_{0}^{\frac{\pi}{2}}} \frac{ 1 \,\,\,\,\,}{(
1-\frac{sin^{2}\theta}{r^2})^{3/2}  }\, d\theta\rightarrow
4\lambda\frac{\pi}{2}=2\pi\lambda$}.
 To prove ($c$.) we compute the first two derivatives of $g(r)$.
A direct calculation shows:

{\footnotesize $$\begin{array}{ccccc} \frac{d}{dr}\,\, g(r)\, =\, 4\lambda\,  r^3
\displaystyle{\int_{0}^{\frac{\pi}{2}}}
\frac{r^2-4sin^2\theta}{\Bigl( r^2-sin^2\theta\Bigr)^{5/2}}\, d\theta, &&&&
\frac{d^2}{dr^2}\,\, g(r)\, =\, 4\lambda \, r^2\,\displaystyle{\int_{0}^{\frac{\pi}{2}}}
\frac{3r^2\, sin^2\theta+12sin^4\theta}{\Bigl( r^2-sin^2\theta\Bigr)^{7/2}}\, d\theta.
\end{array} $$}

It follows that $\frac{d^2}{dr^2}\, g(r)>0 $. This together with
($a$.) and ($b$.) imply that $g(r)$ has a unique minimum at some
point $r_0>1$. Finally we show that $r_0<2$. For this just calculate
{\footnotesize $\frac{d}{dr}\,\, g(2)\, =\, 4\lambda\,  8
{\int_{0}^{\frac{\pi}{2}}} \frac{4-4sin^2\theta}{(
4-sin^2\theta)^{5/2}}\, d\theta\, =\, 32\lambda\,
{\int_{0}^{\frac{\pi}{2}}} \frac{4cos^2\theta}{(
4-sin^2\theta)^{5/2}}\, d\theta\,>0$.} \CaixaPreta \vspace{.3in}

Note that, by definition of $g$, $r_0$ does not depend on $\lambda$
even though $g(r)$ does.
From the Lemma above we can have an idea how the graph of $g(r)$ looks like:\\

\begin{figure}[!htb]
 \centering
 \includegraphics[scale=0.5]{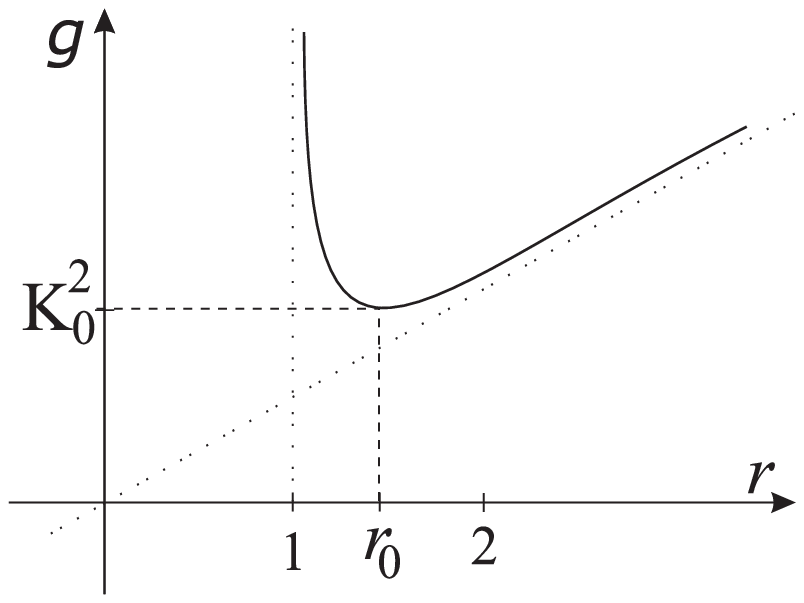}
 \end{figure}

\vspace{.2in}

Define $K_0=\sqrt{g(r_0)}=\sqrt{r_0^3\frac{d}{dr}\, V(r_0)}$.
From the Lemma follows that the restriction of $g(r)$ to the interval
$(1,r_0]$ is a decreasing function. Also $g(r)$ restricted to the interval
$[r_0,\infty )$ is an increasing function. Hence the same is true
for the function $\sqrt{g(r)}$. Therefore the functions
$\sqrt{g}\, |_{(1,r_0]}:(1,r_0]\rightarrow [K_0,\infty )$ and
$\sqrt{g}\, | _{[r_0,\infty )}:[r_0,\infty )\rightarrow [K_0,\infty )$
have inverses which we call $r_1$ and $r_2$ respectively.
It is not difficult to verify that ($iv$), ($v$), ($vii$), ($viii$), ($ix$) hold
for $r_1$ and $r_2$. Next we prove ($vi$). Fix $K>K_0$.
Note that {\small
$U(r)=\frac{K^2}{2r^2}+ V(r)=\frac{K^2}{2r^2}
-4\lambda{\int_{0}^{\frac{\pi}{2}}}
\frac{d\theta} {\sqrt{r^{2}-sin^{2}\theta}}=\frac{1}{r}
\Bigl(\frac{K^2}{2r}
-4\lambda{\int_{0}^{\frac{\pi}{2}}}
\frac{d\theta} {\sqrt{1-\frac{sin^{2}\theta}{r^2}}} \Bigr)$}.
Since $\Bigl( \frac{K^2}{2r}
-4\lambda{\int_{0}^{\frac{\pi}{2}}}
\frac{d\theta} {\sqrt{1-\frac{sin^{2}\theta}{r^2}}} \Bigr)$ tends to $-2\pi\lambda$
when $r\rightarrow\infty$, then for large $r$ we have that $U(r)$ is negative.
This together with ($ii$), ($iii$) of the statement of Theorem C and the fact that
$U(r)$ has exactly two critical points at $r_1$ and $r_2$ imply that
$\frac{d}{dr}\,U(r)>0$ for $1<r<r_1$ and $r>r_2$. To prove that
$\frac{d}{dr}\,U(r)<0$ for $r_1<r<r_2$ it is enough to prove that
$\frac{d}{dr}\,U(r_0)<0$. But
{ $\frac{d}{dr}\,U(r_0)=-\frac{K^2}{r_0^3}+\frac{d}{dr}\, V(r_0)=
-\frac{K^2}{r_0^3}+\frac{K_0^2}{r_0^3}=-\frac{K^2-K_0^2}{r_0^3}<0$}.

Finally we prove ($x$). Let ${\bar r}$ denote $r_1$ or $r_2$.
Since $K^2={\bar r}^3\frac{dV}{dr}({\bar r})$ we have
{\small$U({\bar r})=\frac{K^2}{2{\bar r}^2}+V({\bar r})=
\frac{1}{2}{\bar r}\frac{dV}{dr}({\bar r})+V({\bar r})=
\frac{1}{2}{\bar r}( 4\lambda\int_0^\frac{\pi}{2} \frac{\bar r}{({\bar r}^2-sin^2\theta)^{3/2}} d\theta)
-4\lambda\int_0^\frac{\pi}{2} \frac{d\theta}{({\bar r}^2-sin^2\theta)^{1/2}}=
2\lambda\int_0^\frac{\pi}{2} \frac{2sin^2\theta -{\bar r}^2}{({\bar r}^2-sin^2\theta)^{3/2}}d\theta
=2\lambda\int_0^\frac{\pi}{2} \frac{2cos^2\theta -{\bar r}^2}{({\bar r}^2-cos^2\theta)^{3/2}}d\theta$},
where the last equality is obtained by a change of variable.

If $K\rightarrow +\infty$, we have (by ($viii$) and ($ix$)) that
$r_1\rightarrow 1^+$ and $r_2\rightarrow +\infty$. Hence  it is
enough to prove that $\lim_{r\rightarrow 1} \int_0^\frac{\pi}{2}
\frac{2cos^2\theta -r^2}{(r^2-cos^2\theta)^{3/2}}d\theta =+\infty$
and $\lim_{r\rightarrow +\infty} \int_0^\frac{\pi}{2}
\frac{2cos^2\theta -r^2}{(r^2-cos^2\theta)^{3/2}}d\theta =0$. The
last limit is clearly zero.
 We prove now
$\lim_{r\rightarrow 1}
\int_0^\frac{\pi}{2} \frac{2cos^2\theta -r^2}{(r^2-cos^2\theta)^{3/2}}d\theta =+\infty$.
We have $\int_0^\frac{\pi}{2} \frac{2cos^2\theta -r^2}{(r^2-cos^2\theta)^{3/2}}d\theta =
\int_0^\frac{\pi}{6} \frac{2cos^2\theta -r^2}{(r^2-cos^2\theta)^{3/2}}d\theta +
\int_\frac{\pi}{6}^\frac{\pi}{2} \frac{2cos^2\theta -r^2}{(r^2-cos^2\theta)^{3/2}}d\theta$.
Since
 $\lim_{r\rightarrow 1}
\int_\frac{\pi}{6}^\frac{\pi}{2} \frac{2cos^2\theta -r^2}{(r^2-cos^2\theta)^{3/2}}d\theta=
\int_\frac{\pi}{6}^\frac{\pi}{2} \frac{2cos^2\theta -1}{sin^3\theta}d\theta$
is finite, it is enough to prove that
$\lim_{r\rightarrow 1} \int_0^\frac{\pi}{6} \frac{2cos^2\theta
-r^2}{(r^2-cos^2\theta)^{3/2}}d\theta=\infty$. We can assume
$r<\sqrt{\frac{5}{4}}$. Hence $2cos^2 \theta -r^2 \geq
\frac{1}{4}$, for $\theta\in[0,\frac{\pi}{6}]$.
 Therefore $
\int_0^\frac{\pi}{6} \frac{2cos^2\theta
-r^2}{(r^2-cos^2\theta)^{3/2}}d\theta\geq
\frac{1}{4}\int_0^\frac{\pi}{6}
\frac{d\theta}{((r^2-1)+sin^2\theta)^{3/2}}$. But since
$(r^2-1)<\frac{1}{4},\,\, sin^2\theta\leq sin\theta <1$ and
$sin\theta\leq \theta$ we have $(r^2-1)+sin^2\theta < \frac{5}{4}$.
Hence ${((r^2-1)+sin^2\theta)^{3/2}}< \frac{\sqrt
5}{2}((r^2-1)+sin\theta)\leq \frac{\sqrt 5}{2}((r^2-1)+\theta)$.

Consequently $\int_0^\frac{\pi}{6}
\frac{d\theta}{((r^2-1)+sin^2\theta)^{3/2}}\geq \frac{2}{\sqrt
5}\int_0^\frac{\pi}{6} \frac{d\theta}{(r^2-1)+\theta}=
\frac{2}{\sqrt5}[ln((r^2-1)+\frac{\pi}{6})-ln(r^2-1)]\rightarrow\infty$
as $r\rightarrow 1$. \CaixaPreta \vspace{.3in}

Before we prove Proposition 3 we need some definitions and comments.
Since $V({\mbox\rv})$, ${\mbox\rv}\in\R^3-\, \cal C$, is invariant by rotations around the
$z$-axis we can reduce our problem in a canonical way to a problem with
two degrees of freedom.
Using  cylindrical coordinates ($r,\varphi,z$), the Lagrangian in
these coordinates can be written as:
$L(r,\varphi,z,\dot r,\dot\varphi,\dot z)=\frac{1}{2}({\dot
r}^{2}+r^{2}\dot\varphi ^{2}+\dot z ^{2})-V(r,z)$,
where
$V(r,z)=V({\mbox\rv})$, $ {\mbox\rv}=(rcos\varphi , r sin\varphi,z)$.
Then it is straightforward to verify that the system
$\stackrel{..}{\mbox{\rv}}\,=-\nabla V({\mbox{\rv}})$
in these coordinates is given by:

{\footnotesize \begin{equation}
\left\{ \begin{array}{l}
\stackrel{..}{r}\,  =\displaystyle\frac{K^{2}}{r^{3}}
-\frac{\partial V}{\partial r}
\left( {r},{z}\right)\\  \\
\stackrel{..}{z}\,  =-\displaystyle\frac{\partial V}{\partial z}
\left( r,z\right)\\ \\
\dot\varphi=\displaystyle\frac{K}{r^{2}}
\end{array}
\right.
\label{R4}
\end{equation}}

\noindent where $K$ is the (constant) angular momentum. Note that
the first two equations of system (\ref{R4}) can be rewritten  as
$(\stackrel{..}{r},\stackrel{..}{z} )\,  =-\nabla {\overline
U}\Bigl(r,z\Bigr),$ with { $\overline U
\Bigl(r,z\Bigr)=\frac{K^{2}}{2 r^{2}} +\,V \left(
{r},0,{z}\right).$} If $(r(t),z(t))$ is a solution of the first two
equations of (\ref{R4}), defining {\footnotesize
$\varphi(t)=\int_{0}^{t}\frac{K \,ds}{r^{2}(s)},$} we have that
$(r(t),\varphi(t),z(t))$ is a solution of (\ref{R4}). Then
$(r(t)\,cos\varphi(t),r(t) \, sin\varphi(t),z(t))$ is a solution of
$\stackrel{..}{\mbox{\rv}} =-\nabla V({\mbox{\rv}})$. Note that if
$z\equiv 0$, system (\ref{R4}) is reduced to (\ref{6}) and
$\overline U$ becomes the effective potential $U$.

\vspace{0,4cm}

\noindent {\bf Proof of Proposition 3.} Suppose ${\mbox{\sv}}=(r,z)$ is an equilibrium position of
$\stackrel{..}{\mbox\sv} \,  =-\nabla {\overline
U}(\mbox\sv )$. Then $\nabla {\overline
U}(\mbox\sv )=0$. In particular $\frac{\partial V}{\partial z}=0$. But
$\frac{\partial V}{\partial z}=\beta
z,$ where $\beta =\int_{{\cal
C}}\frac{\lambda\,\,du}{|\!|{\mbox{\rv}}-u|\!|^3}>0$. Hence $z=0$, i.e the solution
lies in the horizontal plane. Since $z=0$, we have that $\overline U (r,0)=U(r)$, where
$U$ is the effective potential.

Let $(r,0)$ be an equilibrium position of $\stackrel{..}{\mbox\sv}
\,  =-\nabla {\overline U}(\mbox\sv )$. If $1<r<r_0$ we know that
the corresponding circular solution is not stable. To prove that
the circular solution that corresponds to $ (r,0)$, $r>r_0$, is
stable it is enough to prove that $(r,0)$ is a strict local
minimum of $ \overline U$.

Since $(r,0)$ is a critical point of $\overline U$, we have that
$r$ is a critical point of $U$, and since $r>r_0$, $r=r_2(|K|)$ where $K$
is the angular momentum of the solution.
 Hence $r$ is a strict local minimum  of $U(r)$ that is $U(r)<U(r')$ for
$r'$ close to $r$, $r'\neq r$. But an
easy  calculation from the definition of $V$ shows that
$V(r',0,0)\leq V(r',0,z)$, for all $r',z$, $(r',z)\not= (\pm 1,0)$
and the same holds for $\overline U$. This
 implies that  $\overline U(r',z)\geq\overline U(r',0)=U(r')>U(r)=\overline U(r,0)$
for $r'$ close to $r$, $r\neq r'$.
Therefore $r=r_2(|K|)$ is a strict local minimum of $\overline U$ .
\CaixaPreta \vspace{0,5cm}

To finish this section we prove that the origin is the only equilibrium position.
Note that by symmetry it is a simple matter to show that the origin is in fact
an equilibrium position.

{\prop  The origin is only equilibrium solution of the  system
$\stackrel{..}{\mbox{{\mbox{\rv}}}}\,  =-\nabla V({\mbox{{\mbox{\rv}}}})$.}
\vspace{.1in}

\noindent  {\bf  Proof.} For $ {\mbox{\rv}}=(x,y,z)$
differentiating we have that $\frac{\partial V}{\partial
z}({\mbox{\rv}})=\beta z,$ where $\beta =\int_{{\cal
C}}\frac{\lambda\,\,du}{|\!|{\mbox{\rv}}-u|\!|^3}>0$. Then if
$\nabla V({\mbox{\rv}})=0$ then $z=0.$ Hence every equilibrium
position lies in the  $xy$-plane. The Proposition now follows from
($i$) of Theorem B and ($i$) of Theorem C. \CaixaPreta \vspace{.4in}

\section{Proof of Proposition 4. }

Fix $\lambda>0$. Consider the  fixed homogeneous
circle ${\cal C}_\epsilon$ with constant density  $\lambda$, contained in the $xy$-plane
and {\it  passing
through the  origin}, with radius $\frac{1}{\epsilon}$ and center $(\frac{1}{\epsilon},0,0)$.
Note that  the $xz$-plane is an invariant subspace of our problem.
Denote  the  potential of this   translated fixed homogeneous circle,
and restricted to  the $xz$-plane,  by
$W\left(x,z;{\epsilon}\right)$. For notational purposes in what follows we use
coordinates $(x,y)$ instead of coordinates $(x,z)$.

 This   potential can be written in the form
$W\left(x,y;\epsilon\right) = -4\lambda\frac{1}{\epsilon}\,T(D,d)$
where $T, D, d$ are as in section 1.1.
If $x\leq\frac{1}{\epsilon}$ we have that
$D^{2}=(\frac{2}{\epsilon}-x)^2 +y^2$ and $d^{2}=x^2 +y^2 $. Hence we have (see section 1.1)
$W\left(x,y;\epsilon\right)= -4\lambda \frac{c}{2} f\left(\frac{x^2
+y^2 }{(c-x)^2 +y^2}\right) / \sqrt{(c-x)^2 +y^2}$, where $c=\frac{2}{\epsilon}$.

Let $\nabla W$ denote the gradient $(\frac{\partial W}{\partial x},\frac{\partial W}{\partial y})$
of $W$. A straightforward calculation shows that for $x\leq c=\frac{2}{\epsilon}$ we have

{\lem
$$\nabla W\left(x,y;\epsilon\right) = 2\lambda c
\frac{f(t)}{D^3}(x-c,y) - 4\lambda \frac{c^2}{D^5}f'(t) \,[c(x,y) +(y^2
-x^2 ,-2xy)],$$
\noindent where $t=\frac{x^2 +y^2}{(c-x)^2 +y^2}$  and $D= \sqrt{(c-x)^2
+y^2}$.
\label{1.6.3}}

\vspace{0,6cm}

\noindent{\bf Remark.} If  $z= x+iy$ then $z^2 =(x^2 -y^2, 2xy)$, hence the identity above becomes
$$\begin{array}{l}\,\,\,\,\,\,\,\,\nabla W\left(z;\frac{c}{2}\right) = 2\lambda
\frac{c}{|z-c|^3}f\left(\frac{|z|^2}{|z-c|^2}\right)(z-c) - 4\lambda
\frac{c^2}{|z-c|^5}f'\left(\frac{|z|^2}{|z-c|^2}\right)(cz-z^2)\,=\\
\\
=-2\lambda
\frac{c^{2}}{|z-c|^3}f\left(\frac{|z|^2}{|z-c|^2}\right) +\left[ 2\lambda
\frac{c}{|z-c|^3}f\left(\frac{|z|^2}{|z-c|^2}\right) - 4\lambda
\frac{c^3}{|z-c|^5}f'\left(\frac{|z|^2}{|z-c|^2}\right)\right] z +4\lambda
\frac{c^2}{|z-c|^5}f'\left(\frac{|z|^2}{|z-c|^2}\right) z^{2}\\\\
=-2\lambda
\frac{c}{|z-c|^3}\Bigl\{
f\left(\frac{|z|^2}{|z-c|^2}\right) +z\left(
\frac{2 c}{|z-c|^2}\right)f'\left(\frac{|z|^2}{|z-c|^2}\right)
\Bigr\}\,(z -c).
\end{array}$$

\vspace{0,15cm}

\noindent Note that the   gradient vector has components in the
$z$ and $z^{2}$ directions and in the
 real $(1,0)$ direction.\\

Now set
$h_{1}(x,y;\epsilon)=\displaystyle\frac{4\lambda}
{\epsilon}\frac{f(t)}{D^{3}} \frac{1}{\epsilon ^{3/2}}$,
$\,h_{2}(x,y;\epsilon)=-\displaystyle\frac{32\lambda}{{\epsilon^{3}}}\frac{1}{D^{5}}
\frac{1}{\epsilon^{2}}\,$ and
$\,h_{3}(x,y;\epsilon)=\displaystyle\frac{16\lambda}{{\epsilon^{2}}}\frac{f'(t)}{D^{5}}
\frac{1}{\sqrt{\epsilon} }$.
Hence
$$\nabla W(x,y;\epsilon )=
[\epsilon ^\frac{3}{2} h_1] (x,y) -[2\sqrt{\epsilon}\,h_1 ] (1,0)
+ [\epsilon ^2 f'(t) h_2] (x,y)+ [\sqrt{\epsilon}\, h_3]  (x^2
-y^2 ,2xy).$$

\vspace{0,3cm}

\noindent The next Lemma shows that $h_1, h_2, h_3$ are bounded
in certain sense.

{\lem
Given a compact $C\subset \R ^2-\{0\}$
 there are $\epsilon_{0}>0$, $K>0$ such that $|h_i(x,y;\epsilon)|<K$
for all $(x,y)\in C$, $0<\epsilon\leq\epsilon_{0}$. \label{6.2}}

\vspace{0,3cm}

\noindent  {\bf  Proof.} Given a compact $C\subset\R^{2}-\{0\}$,
there exist constants $C_{1},C_{2},$ such that $0<C_{1}\leq
\|(x,y)\|=\sqrt{x^{2}+y^{2}}\leq C_{2},$ for all $(x,y)\in C$. For
$\epsilon_{0}<\frac{1}{\mbox{$C_{2}$}},$ we have
$C_{2}<\frac{1}{\mbox{$\epsilon_{0}$}},$ which implies that
$C\subset \{ (x,y);x< \frac{1}{\epsilon}, y\in\R\}$, for all
$\epsilon\leq \epsilon_{0},$
 and $h_{i}(x,y;\epsilon)$ is defined in $C$, for
$\epsilon\leq \epsilon_{0},\,\,i=1,2,3.$

Since $\sqrt{(2-\epsilon x)^{2}+\epsilon^{2}y^{2}} =
\|(2,0)-\epsilon \,(x,y)\|\,\leq \,
\|(2,0)\|+\epsilon\,\|(x,y)\|\,= \,2+\epsilon \,\|(x,y)\|\,$ and
$\|(2,0)-\epsilon \,(x,y)\|\, \geq \, \|(2,0)\|-\epsilon\,\|(x,y)\|\,=\,2-
\epsilon
\,\|(x,y)\|,$ and  we obtain

\begin{equation}
 2-\epsilon\, C_{2}\leq  \sqrt{(2-\epsilon x)^{2}+\epsilon^{2}y^{2}}
\leq 2+\epsilon\, C_{2}
\label{5.3}
\end{equation}

Therefore, using  Lemma \ref{1.6.2}, part (iii), and
(\ref{5.3}) we have that:

{\footnotesize{$$|h_{1}(x,y;\epsilon)|={\left|
\frac{4\lambda}{\epsilon}\frac{f(t)}{D^{3}}\frac{1}{\epsilon ^{3/2}}\right|}
= \frac{4\lambda\sqrt{\epsilon}\,f(t)}{((2-\epsilon
x)^{2}+\epsilon^{2}y^{2})^{3/2}}
\leq\frac{4\lambda\sqrt{\epsilon}\,f(t)}{(2-\epsilon C_{2})^{3}}
 \,\leq\, \frac{2\pi \lambda\sqrt{\epsilon}}{(2-\epsilon C_{2})^{3}\sqrt[4]{t}}=$$}}
{\footnotesize{$$ = \displaystyle\frac{2\pi \lambda\sqrt{\epsilon }\,[(2-\epsilon
x)^{2}+\epsilon^{2}y^{2}]^{1/4}}{(2-\epsilon C_{2})^{3}\sqrt{\epsilon}{(x^{2}+y^{2})^{1/4}}}
\leq   \displaystyle\frac{2\pi\lambda (2+\epsilon C_{2})^{1/2}}{( 2-\epsilon
C_{2})^{3}C_{1}^{1/2}}
\leq \displaystyle\frac{2\pi\lambda (2+\epsilon_{0} C_{2})^{1/2}}{(2-\epsilon_{0}
C_{2})^{3}C_{1}^{1/2}}=K_{1}(C,\epsilon_{0}).$$}}

This proves the Lemma for $h_{1}$.

Now, applying (\ref{5.3}) we have
{\small ${|h_{2}(x,y;\epsilon)| =
\frac{32\lambda}{D^{5}\epsilon^{5}}
= \frac{32\lambda}{((2-\epsilon x)^{2}+\epsilon^{2}y^{2})^{5/2}}
\leq \frac{32\lambda}{(2-\epsilon C_{2})^{5}}\leq\frac{32\lambda}
{(2-\epsilon_{0}C_{2})^{5}}=}$}
{\small $K_{2}(C,\epsilon_{0}),$}
 and this proves the Lemma for $h_{2}$.

Finally, using  Lemma \ref{1.6.2}, part (v) and (\ref{5.3}), we have
 $|h_{3}(x,y;\epsilon)| =\left|
\frac{16\lambda}{{\epsilon^{2}}}\frac{f'(t)}{D^{5}}
\frac{1}{\sqrt{\epsilon} }\right|$
$=\frac{16\lambda}{\epsilon^{5/2}}\frac{|f'(t)|}{\left[
\frac{(2-\epsilon x)^{2}+\epsilon^{2}y^{2}}{\epsilon^{2}}\right]^{5/2}}$
 $= \frac{16\lambda \epsilon^{5/2}|f'(t)|}{((2-\epsilon x)^{2}+
\epsilon^{2}y^{2})^{5/2}}
\leq \frac{4\lambda\pi}{(x^{2}+y^{2})^{5/4}}\frac{1}{((2-\epsilon x)^{2}+
\epsilon^{2}y^{2})^{5/4}}
\leq \frac{4\lambda\pi}{C_{1}^{5/2}}\frac{1}{(2-\epsilon
C_{2})^{5/2}}
\leq  \frac{4\lambda\pi}{C_{1}^{5/2}}\frac{1}{(2-\epsilon_{0}
C_{2})^{5/2}}= K_{3}(C,\epsilon_{0}).$ \CaixaPreta\\

\vspace{0,6cm}

Write $h(x,y;\epsilon)=(x^2 +y^2) \epsilon^2 f'(t) h_2$.
Then $\nabla W= [\epsilon ^\frac{3}{2} h_1] (x,y)
-[2\sqrt{\epsilon}\,h_1 ] (1,0) + h \,\frac{(x,y)}{(x^2 +y^2)}+
[\sqrt{\epsilon}\, h_3]  (x^2 -y^2 ,2xy)$.

{\lem If $\,\,(x_n,y_n,\epsilon_n)\rightarrow (x,y,0)\neq
(0,0,0),$ with  $\epsilon_n\neq 0,\,\,(x_n,y_n)\notin {\cal
C}_{\epsilon_{n}},$ then \\$lim_{n\rightarrow +\infty}
h(x_n,y_n;\epsilon_n)=64\lambda.$} \vspace{0,3cm}

\noindent {\bf Proof.} First note that $\lim_{n\rightarrow
+\infty} h_2(x_n,y_n;\epsilon_n)=-32\lambda$. Also
$\lim_{n\rightarrow +\infty} \frac{\epsilon_n ^2}{t_n}
=\frac{4}{x^{2}+y^{2}}$, where
 $ t_n = \frac{x_n^2  +y_n^2}{\left(x_n -\frac{2}{\epsilon_n}\right)^2
+y_n^2}.\,\,$ Therefore  $\lim_{n\rightarrow +\infty}
h(x_n,y_n,\epsilon_n)=\lim_{n\rightarrow +\infty} (x_n ^2 +y_n
^2)(\frac{\epsilon_n ^2}{t_n})(f'(t_n)
t_n)h_2(x_n,y_n;\epsilon_n)=(x^2 +y^2)
\frac{4}{x^{2}+y^{2}}(-\frac{1}{2})(-32\lambda)=64\lambda$,
where we are using the fact that $\lim_{t\rightarrow 0^{+}}
tf(t)=-\frac{1}{2}$ (see Lemma \ref{1.6.2}).
 \CaixaPreta

\vspace{0,6cm}

 Define $\nabla W(x,y;0)
:=64\,\lambda\frac{(x,y)}{x^{2}+y^{2}},$ $(x,y)\neq (0,0).$ Hence
$\nabla W(x,y;\epsilon)$ is defined on $A=\{\,(x,y;\epsilon)\,;\,
x^{2}+y^{2}\neq 0,\,(x,y)\neq (\frac{2}{\epsilon},0)\, \}$, and
for each  $\epsilon$, $\nabla W$ is analytic.

\vspace{0,4cm}

\noindent {\bf Proof of Proposition 4.} Clearly $\nabla W$ is
continuous in $(x,y,\epsilon)$, with $\epsilon\neq 0$. Let
$\{\,(x_{n}, y_{n}, \epsilon_{n})\,\}$ be a sequence in $A$, with
$(x_{n}, y_{n}, \epsilon_{n})\rightarrow (x,y,0),$ $\,\,(x,y)\neq
(0,0),\,\,\epsilon_n\neq 0.$ Thus $C=\{\,(x_{n},
y_{n}),\,n\in\N\,\}\,\cup\,\{\,(x,y)\,\}$ is compact and
$(0,0)\notin C$. By Lemma \ref{6.2} there exist $n_{0}$ and
$K_{C}>0$, such  that
  $\,|h_i((x_{n}, y_{n},
\epsilon_{n})|<K_{C},$ for $n\geq n_{0},$ $i=1,2,3.$
Therefore, $ \lim_{n\rightarrow +\infty}\nabla W (x_{n},
y_{n},\epsilon_{n})=\lim_{n\rightarrow +\infty} h(x_n,y_n;\epsilon_n)
\frac{(x_n,y_n)}{x_n^2+y_n^2}=64\lambda \frac{(x,y)}{x^2 +y^2}$.
 \CaixaPreta

\vspace{0,5cm}

\section{ Proof of Proposition 5.}

Take $M=1$. Let
$\beta(t)=(D(t),d(t))$ with $\beta(0)=(D,d)$, $\dot{\beta}(0)=(u,v
).$ Let us calculate $\sigma(t)=\sigma(D(t),d(t))$. The following
change of variables will simplify our calculations. Write
$u=\gamma D$ and $v=(\gamma +\delta)d$. The Taylor series of
$D(t),d(t)$ around $t=0$ are:

$$D(t)= D+tu +\ldots =(1+t\gamma )D +O(t^{2})$$
$$d(t) =d+tv +\ldots =(1+t\gamma )d + t\delta d + O(t^{2})$$

\vspace{0,3cm}

{\lem  $\sigma (t)=(1+t\alpha^\ast)\sigma + O(t^2),$ with
$\alpha^\ast =\alpha_0 +\delta\sum_{n=1}^{+\infty} \frac{1}{2^{n}}
      \frac{ d_{n-1}}{D_{n}}\prod_{j=1}^{n-1}\left(
      1-\frac{d_{j-1}}{D_{j}}\right)$,
 $\alpha _{0}=\gamma$,\,\, $d_{0}=d$,\,\, $D_{0}=D$, $D_{n}=
\frac{D_{n-1}(0)+d_{n-1}(0)}{2}$, $d_{n}=\sqrt
{D_{n-1}(0)d_{n-1}(0)}$, for $n\geq 1$.}

\vspace{0,3cm}

\noindent {\bf Proof.} Calculating the arithmetic-geometric mean
of $D(t)$ and $d(t)$, we have:
$$D_{1}(t) = \frac{D(t)+d(t)}{2}
         = (1+t\gamma )D_{1} +\frac{1}{2} t\delta d + O(t^{2})
         = (1+t\alpha_{1} )D_{1} + O(t^{2}),$$
$$d_{1}(t)=\sqrt{D(t)d(t)} =(1+t\gamma )d_{1} +\frac{1}{2} t\delta
d_{1} + O(t^{2}) =  (1+t\alpha _{1} )d_{1} +\frac{1}{2} t\delta
       d_{1}\left( 1-\frac{d}{D_{1}}\right) + O(t^{2}),$$
\noindent where $\alpha_{1}=\gamma +\frac{1}{2}\delta
\frac{d}{D_{1}}$,  $D_{1} = \frac{D(0)+d(0)}{2} = \frac{D+d}{2}\,\,$ and $\,d_{1}=\sqrt
{D(0)d(0)}=\sqrt{Dd}$.\\

Calculating the arithmetic and geometric mean successively,
rearranging and defining\\
 $\alpha _{0}=\gamma$, \,\,$\alpha_{1}=\gamma
+\frac{1}{2}\delta \frac{d}{D_{1}}$ and
$$\alpha_{n}= \alpha_{n-1} + \frac{1}{2^{n}} \delta
      \frac{ d_{n-1}}{D_{n}}\prod_{j=1}^{n-1}\left(
      1-\frac{d_{j-1}}{D_{j}}\right),\,\,\mbox{for}\,\,\,\,
      n\geq 2 $$
\noindent we obtain, by induction, the general rule:
$$D_{n}(t) = (1+t\alpha_{n} )D_{n} +O(t^{2}),$$
$$d_{n}(t)=  (1+t\alpha _{n} )d_{n} +\frac{1}{2^{n}} t\delta
       d_{n}\prod_{j=1}^{n}\left( 1-\frac{d_{j-1}}{D_{j}}\right) +
       O(t^{2}),
$$
\noindent where  $d_{0}=d$,\,\, $D_{0}=D$, $D_{n}=
\frac{D_{n-1}(0)+d_{n-1}(0)}{2}$, $d_{n}=\sqrt
{D_{n-1}(0)d_{n-1}(0)}$, \,\,$n\geq 1$.

\noindent Note that $ \alpha_{n} =\alpha_0 +\delta\sum_{i=1}^{n}
 \frac{1}{2^{i}}  \frac{ d_{i-1}}{D_{i}}\prod_{j=1}^{i-1}\left(
      1-\frac{d_{j-1}}{D_{j}}\right)$.

Finally, taking the limit on $d_{n}(t)$ or on $D_{n}(t)$, when $n$
goes to
      infinity, we have $$\sigma (t)=(1+t\alpha^\ast)\sigma + O(t^2),$$
with $\alpha^\ast =\lim_{n\rightarrow\infty}\alpha_n=\alpha_0
+\delta\sum_{n=1}^{+\infty} \frac{1}{2^{n}}
      \frac{ d_{n-1}}{D_{n}}\prod_{j=1}^{n-1}\left(
      1-\frac{d_{j-1}}{D_{j}}\right)$.
It is easy to see that this series converges. \CaixaPreta\\

\vspace{0,6cm}

\noindent{\bf Proposition 5.} {\it We have the following formulas:
$$\frac{\partial}{\partial D}V(D,d)=\frac{\chi -1}{D}
V(D,d),\,\,\,\,\,\,\,\,\,\, \frac{\partial}{\partial d}
V(D,d)=\,-\frac{\chi}{d}V(D,d),$$

\noindent where $\chi=\sum_{n=1}^{+\infty} \frac{1}{2^{n}}
      \frac{ d_{n-1}}{D_{n}}\prod_{j=1}^{n-1}\left(
      1-\frac{d_{j-1}}{D_{j}}\right)$.}\\

\noindent {\bf Proof.} Since $V(P(t))= \frac{1}{\sigma(t)}=
\frac{1}{\sigma} -\frac{t\alpha^{\ast}}{\sigma}+O(t^{2}) = V(D,d)
-t\alpha^{\ast} V(D,d) +O (t^{2})$  we can now calculate the partial
derivaties $\frac{\partial}{\partial d}V(D,d)$,
$\frac{\partial}{\partial D}V(D,d)$.

We have $\frac{d}{dt} V(\beta(t))_{|t=0} = -\alpha^{\ast} V(D,d) =
\langle \,(\frac{\partial}{\partial D}V(\beta(0)),
\frac{\partial}{\partial d}V(\beta(0))), \,(u,v)\,\rangle$. Note that
$\alpha^{\ast}=\gamma +\delta \chi$, with
$\chi
=\sum_{n=1}^{+\infty} \frac{1}{2^{n}}
      \frac{ d_{n-1}}{D_{n}}\prod_{j=1}^{n-1}\left(
      1-\frac{d_{j-1}}{D_{j}}\right)$, $\gamma =\frac{u}{D}$ and
$\,\delta =\frac{v}{d}-\frac{u}{D}$. We have $\alpha^{\ast}
=\,\frac{u}{D} +\left( \frac{v}{d}-\frac{u}{D}\right)\chi\,
=\,\frac{u}{D}(1-\chi) +\frac{v}{d}\chi.$
 Therefore
 $\alpha^{\ast} =\left(\frac{1-\chi}{D}\right) u
+\frac{\chi}{d}v, $ so $\langle\, (\frac{\partial}{\partial
D}V(D,d), \frac{\partial}{\partial d}V(D,d)), \,(u,v)\,\rangle =-
\left(\frac{1-\chi}{D} V(D,d)u+ \,\frac{\chi}{d}V(D,d)v\right)$ for all $u, v$.
Hence:
$$(\frac{\partial}{\partial
D}V(D,d), \frac{\partial}{\partial d}V(D,d))=\left( \frac{\chi
-1}{D} V(D,d), \,-\frac{\chi}{d}V(D,d)\right).\, \CaixaPreta
$$


{}


\vspace{.3in}



\begin{thebibliography}{xxx}

\bibitem[1]{AO} C. Azev\^edo and P. Ontaneda, {\em On the existence of
periodic orbits for the
 fixed homogeneous circle problem}.
Submitted for publication. ArXiv: math.CA/0502305.

\bibitem[2]{BLO} E. Belbruno,   J. Llibre and  M. Oll\'e,
{\em `On the families of periodic orbits which bifurcate from the
circular Sitnikov motions'}, Celestial Mechanics and  Dynamical
Astronomy, {\bf 60}, 99-129, 1994.

\bibitem[3] {G}  H. Goldstein, {\em Classical Mechanics,}
Addison-Wesley Publishing Company, 1980.


\bibitem[4] {Po} H. Poincar\'e, {\em Th\'eorie du Potentiel
Newtonien,} \'Editions Jacques Gabay, Paris,  1990.

\bibitem[5] {S}  J. Sotomayor,  {\em
Li\coes de Equa\coes Diferenciais Ordin\'arias}, Projeto Euclides.

\bibitem[6] {M}  J. Moser, {\em Stable and random motions in dynamical
systems}, Annals of Mathematics Studies  77, Princeton University
Press,  1973.

\bibitem[7] {Mc} W.D. MacMillan, {\em Theoretical
Mechanics: The theory of the potential,} Dover Publications, Inc.
New York,  1958.

\end {thebibliography}
\end{document}